\documentclass[english]{article}
\usepackage[T1]{fontenc}
\usepackage[latin9]{inputenc}
\usepackage[letterpaper]{geometry}
\usepackage{float}
\usepackage{amsmath}
\usepackage{graphicx}
\usepackage{color}
\usepackage{epstopdf}
\usepackage{multirow}
\usepackage{array}
\usepackage{mathrsfs}
\usepackage{amsmath, amssymb, amsthm}
\usepackage{subfigure}
\usepackage{tikz}

\graphicspath{{./Figures/}}

\newtheorem{remark}{Remark}

\usepackage{amsfonts} 
\newcommand{\Oh}{\mathcal{T}_h}
\newcommand{\Eh}{\mathcal{E}_h}
\newcommand{\EH}{\mathcal{E}_H}

\usepackage{color}
\definecolor{black}{rgb}{0,0,0}

\definecolor{red}{rgb}{1,0,0}

\definecolor{blue}{rgb}{0,0,1}

\newcommand{\specialcell}[2][c]{\begin{tabular}[#1]{@{}c@{}}#2\end{tabular}}
\usepackage{multirow}

\makeatother

\usepackage{babel}

\title{ An enriched multiscale mortar space for high contrast flow problems
}
\author{Eric T. Chung\thanks{Department of Mathematics, The Chinese University of Hong Kong, Hong Kong SAR. } \and
	Shubin Fu\thanks{Department of Mathematics, Texas A\&M University, College Station, TX 77843.}
	\and Yanfang Yang\thanks{Department of Mathematics, Texas A\&M University, College Station, TX 77843.}
}

\begin{document}
	
	\maketitle

\begin{abstract}
Mortar methods are widely used techniques for discretizations of partial differential equations
and preconditioners for the algebraic systems resulting from the discretizations. 
For problems with high contrast and multiple scales, the standard mortar spaces are not robust,
and some enrichments are necessary in order to obtain an efficient and robust mortar space. 
In this paper, we consider a class of flow problems in high contrast heterogeneous media, and 
develop a systematic approach to obtain an enriched multiscale mortar space. 
Our approach is based on the constructions of local multiscale basis functions. 
The multiscale basis functions are constructed from local problems by following the framework of the Generalized Multiscale Finite 
Element Method (GMsFEM). In particular, 
we first create a local snapshot space. Then we select the dominated modes within the snapshot space using
an appropriate Proper Orthogonal Decomposition (POD) technique. These multiscale 
basis functions show better accuracy than polynomial basis for multiscale problems.
Using the proposed multiscale mortar space, we will construct a multiscale finite element method
to solve the flow problem on a coarse grid and a preconditioning technique
for the fine scale discretization of the flow problem. 
In particular, we develop a multiscale mortar mixed finite element method
using the mortar space. 
In addition, we will design a two-level {\it additive} preconditioner and a two-level {\it hybrid} preconditioner based on the proposed mortar space
for the iterative method applied to the fine scale discretization of the flow problem. 
We present several numerical examples to demonstrate the efficiency and robustness of our proposed mortar space with respect to 
both the coarse multiscale solver and the preconditioners.
\end{abstract}


\section{Introduction}
In this paper, we consider the following second order elliptic
differential equation in mixed form:
\begin{subequations}
\begin{alignat}{2} 
\label{original equation-1}
 \boldsymbol{q} + \kappa\nabla u &= 0 \qquad && \text{in $\Omega$,}\\
\label{original equation-2}
\nabla \cdot \boldsymbol{q}  &= f && \text{in $\Omega$,}\\
\label{boundary condition}
u &= 0 && \text{on $\partial \Omega$.}
\end{alignat}
\end{subequations}
where $\Omega \subset \mathbb{R}^d$ ($d=2,3$) is a bounded polyhedral domain with outward unit normal vector $\boldsymbol{n}$ on the boundary, 
$f \in L^2(\Omega)$, $\kappa$ represents the permeability  field that varies over multiple spacial scales. 
Possible applications of (\ref{original equation-1})-(\ref{boundary condition}) include  flows in porous media,  diffusion and transport of passive chemicals or heat transfer
in heterogeneous media.
Solving (\ref{original equation-1})-(\ref{boundary condition}) can be challenging 
if $\Omega$ is large and the 
permeability $\kappa$ is heterogeneous with multiple scales and high contrast, which is a  common characteristic in many  industrial, scientific, engineering, and
environmental applications. Direct simulation requires
very fine meshes and this makes the corresponding algebraic system very large and ill conditioned (due
to both the small mesh size and the high contrast of the coefficient). Thus direct simulation is computationally intractable. 

In order to solve (\ref{original equation-1})-(\ref{boundary condition}) efficiently, various reduced-order methods have been proposed and applied. These methods
include numerical upscaling (see, e.g.,\cite{weh02,dur91}), variational multiscale method (see, e.g.,\cite{hughes98,ref:Iliev2011}),
 multiscale finite element method (see, e.g.,\cite{egw10,efendiev2009multiscale,Arbogast_two_scale_04,chung2013sub}),
  mixed multiscale finite element methods (see, e.g., \cite{chung2015mixed,chen2003mixed}), the multiscale
finite volume method (see, e.g.,\cite{jennylt03}), mortar multiscale finite element method (see, e.g.,\cite{Wheeler_mortar_MS_12,Wheeler_MS_mortar_11, Arbogast_PWY_07,arbogast2013ms}),
multiscale hybrid-mixed finite element methods (see, e.g.,\cite{hybridmixed2013,hybridmixed2015}), generalized multiscale finite element methods
(see, e.g., \cite{efendiev2013generalized,chung2016adaptive,chung2016mixed,gao2016application,chung2015perforated}) 
and weak Galerkin generalized multiscale finite element method \cite{mu2016weak}. 
These methods typically use some type of global couplings in the coarse grid level to link
the sub-grid variations of neighboring coarse regions. 
We will, in this paper, consider the global coupling via the mortar framework. 
The mortar framework offers many advantages, such as the flexibility in the constructions
of the coarse grid and sub-grid capturing tools. The framework also gives a smaller dimensional global system since the degrees of freedom
are reduced to coarse region boundaries. 
The connectivity of the sub-grid variations is typically enforced using a Lagrange multiplier. 
For multiscale problems, the choice of the mortar space for the Lagrange multiplier
requires a very careful construction, in order to obtain an efficient and robust method.
To construct an accurate mortar space with a small dimension, we will apply the recently developed GMsFEM, which offers a systematic approach for model reduction. 
In particular, 
we first create a local snapshot space for every coarse edge. 
We obtain this space by first solving some local problems on a small region containing an edge,
and then restricting the solutions to the edge. 
Next, we select the dominated modes within the snapshot space using
an appropriate POD technique.
These dominated modes form the basis for the mortar space. 
We will apply our mortar space in two related formulations.
The first one is a mixed GMsFEM using the mortar formulation. This method gives
a coarse-grid solver for the problem (\ref{original equation-1})-(\ref{boundary condition}).
The second one is a coarse space for some preconditioners applied to the fine scale discretization of (\ref{original equation-1})-(\ref{boundary condition}).


In multiscale finite element method, the local basis functions are constructed independently in each coarse cell (subregion). In general, the basis functions are discontinuous across the cell interfaces. These discontinuities can be coupled by mortar spaces. The pioneering work on mixed finite element approximations on multiblock grids were introduced and studied by Arbogast et al. in (\cite{arbogast2000mixed}). In the mortar framework, a Lagrange multiplier is introduced to impose the continuity of flux across block interfaces. This Lagrange multiplier lives in an auxiliary space (called mortar space). In a two-scale  method the aim is to resolve the local heterogeneities on the fine grid introduced on each coarse block and then "glue"
these approximations together via mortar spaces.  The global problem is then formulated in terms of the Lagrange multiplier, which in general
yields  a much smaller algebraic system compared with the original systems.

 If the mortar space is small, we can solve the final algebraic system relatively efficiently. However the continuity of the flux is only weakly imposed. Therefore, for both efficiency and accuracy, it is preferable to construct a mortar space that is not too large and can still impose the continuity accurately. Our goal of this work is to construct such a mortar space by a systematic enrichment technique using GMsFEM. We will first create local snapshot spaces and then perform POD to the snapshot spaces to identify important modes to form the enriched multiscale mortar space. In this way, only a few degrees of freedom per interface are used and the basis functions in the space demonstrate a good "gluing" behavior.
 We will also study the effects of using oversampling techniques, randomized snapshots and different
 sizes of local problem domain on the robustness and accuracy of our mortar space.
 Our work share some similarities with multiscale hybridizable discontinuous Galerkin method \cite{ms_hdg1, ms_hdg2}, but the key part is different, which involves the use of a new methodology to construct basis. The mortar mixed finite element enjoys the 
 advantage of global mass conservation, which is important in industrial reservoir simulations.

On the other hand, based on our proposed mortar space,
we will construct an effective and robust two-level preconditioner to solve the algebraic system arising from a 
fine scale discretization by some iterative methods. 
Our approach uses the solutions of small local problems and a coarse problem in constructing the preconditioners for the fine-scale system. It is well known that
for high contrast heterogeneous media,
if the coarse problem is not suitably chosen, the performance of the preconditioner may deteriorate. To deal with this problem,  many researchers designed different types of robust two-level preconditioners with nonstandard coarse problems in the past
several decades. For example, in \cite{sarkis1997dd} the authors proposed a nonstandard coarse space for the elliptic problems with discontinuous coefficients.
The idea of using single multiscale basis to form the coarse space is reported in
\cite{graham2007dd,oh2013overlapping}, this method is robust if the high conductivity
does not cross the coarse grid. Using spectral functions to enrich the coarse space
turns out to be very efficient and robust for problems with almost any types of media
(see, e.g.,\cite{galvis2010domain,galvis2010domain2,efendiev2012robust,dolean2012analysis,kim2015bddc,spillane2014abstract,kim2016bddc,kim2016sdg,klawonn2015feti,oh2016bddc}). 

Following previously mentioned works,
we propose a two-level preconditioner which utilizes the multiscale mortar space designed in this paper to construct the coarse preconditioner. The local preconditioner is formed 
by solving a local Dirichlet problem in the neighborhood of each edge.
This local preconditioner will be block diagonal and is more effective for highly heterogeneous coefficient than the coarse element based local preconditioner, see \cite{xiao2013multiscale}. However, neighborhood defined local preconditioner may
lead to higher computational cost. To solve this problem, we will also consider
restrictive local preconditioner \cite{cai1999restricted}, which is quite similar
to the oversampling techniques when constructing the multiscale basis. 
Moreover, we will study the effects of the size of local problems and robustness of the 
method.
We remark that the preconditioners designed in this paper is similar to the one in \cite{arbogast2015two, xiao2013multiscale}, the major difference
is the construction of the coarse space. In \cite{arbogast2015two}, polynomials and homogenized basis functions are applied 
for the coarse preconditioner, which shows good performance for checker board type high-contrast model, however 
there is no evidence that the polynomial or homogenized basis functions can deal with arbitrary type of
high contrast media.
 Our numerical results show that the proposed coarse space gives promising ability to deal with problems in
 media with complicated inclusions and long channels that may across coarse edges.
 
The paper is organized as follows. In section \ref{Form4meth}, we first describe the coarse and fine discretizations of the domain, then present the framework of mortar mixed finite element method, followed by the description of the domain decomposition method. In section \ref{ms-basis}, we introduce the construction  of the multiscale basis. In section \ref{two-level}, the coarse and local preconditioners are defined, which are combined to form our two-level preconditioners. Numerical examples are given in section \ref{numerical-results}, and conclusions are made in the last section.

\section{Preliminaries} \label{Form4meth}

In this section, we will give some basic definitions. 
We will also present the formulations of a fine-scale discretization for (\ref{original equation-1})-(\ref{boundary condition})
and its domain decomposition formulation, as well as the formulation for a mixed multiscale method for (\ref{original equation-1})-(\ref{boundary condition}).

\subsection{Fine and coarse grids}
The proposed multiscale mortar space requires a coarse grid and a fine grid, which will be introduced in this section.
 We first divide the computational domain $\Omega$ into non-overlapping polygonal coarse blocks  
 $K_i$ with diameter $H_i$ so that $\overline{\Omega}=\cup_{i=1}^N\overline{K}_i$,
 where $N$ is the number of coarse blocks.
The decomposition of the domain can be nonconforming. We call $E$ a coarse
 edge of the coarse block $K_i$ if $E = \partial K_i \cap \partial K_j $ or $E_H= \partial K_i \cap \partial{\Omega}$.
 Let $\mathcal{E}_H(K_i)$ be the set of all coarse edges on the boundary of the coarse block $K_i$ and $\mathcal{E}_H=\cup_{i=1}^N\mathcal{E}_H(K_i)$
 be the set of all coarse edges.

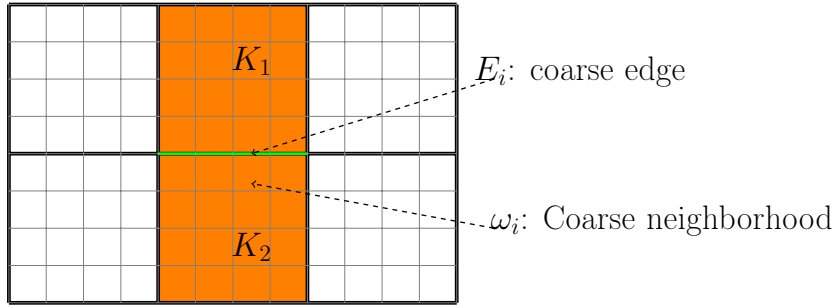
\begin{figure}[H]
\centering
\resizebox{0.7\textwidth}{!}{
\begin{tikzpicture}[scale=0.8]
\filldraw[fill=orange, draw=black] (4,0) rectangle (8,8);
\draw[step=4,black, line width=0.8mm] (0,0) grid (12,8);
\draw[step=1,gray, thin] (0,0) grid (12,8);
\draw[ultra thick, green](4, 4) -- (8,4);
\draw [->,dashed, thick](13,6) -- (6.5, 4);
\node at (15.3,6.2)  {\huge $E_{i}$: coarse edge};

\draw [->,dashed, thick](13,2) -- (6.5, 3.2);
\node at (17.5,2.2)  {\huge $\omega_{i}$: Coarse neighborhood };
\node at (6.5,6.5)  { \huge $K_{1}$} ;
\node at (6.5,1.5)  { \huge $K_{2}$} ;
\end{tikzpicture}
}
\caption{ Illustration of the coarse and fine grids and 
a coarse neighborhood.}
\label{grids}
\end{figure}

For each coarse block $K_i$, we introduce a shape regular discretization $\mathcal{T}_h(K_i)$  with rectangular elements (denoted as $T_j, j=1,\cdots, i_n$) with mesh size $h_i$. 
Let $\mathcal{T}_h=\cup_{i=1}^N\mathcal{T}_h(K_i)$ be the union of all these triangulations, which is a fine mesh triangulation of the domain $\Omega$.
We use $h$ to denote the mesh size of $\mathcal{T}_h$. 
In addition, we let
$\mathcal{E}_h(K_i)$
be the set of all edges of the triangulation $\mathcal{T}_h(K_i)$ and $\mathcal{E}_h^0(K_i)$ be the set of all interior edges of the triangulation
$\mathcal{T}_h(K_i)$ and let $\mathcal{E}_h=\cup_{i=1}^N\mathcal{E}_h(K_i)$ be the set of all edges in the triangulation $\mathcal{T}_h$.
We denote the edges of this triangulation by generically by $e$. 
Figure~\ref{grids} gives an illustration of the constructions of the two grids. The black lines represent the coarse grid, and the grey lines represent the fine grid. 
For each coarse edge $E_i$, we define a coarse neighborhood $\omega_i$ as the union of all coarse blocks having the edge $E_i$.
Figure~\ref{grids} shows a coarse neighborhood $\omega_i$ in orange color. 

\subsection{Fine scale mortar mixed finite element method}
In this section, we recall the standard mortar mixed finite element method defined on the fine grid. 
First of all, 
using the two-scale grid defined in the last subsection, we define the following finite element spaces:
\begin{alignat*}{1}
W_h:=&\;\{w\in{L}^2(\Oh):\;w|_T\in P^0(T),\; T\in\Oh\}, 
\\
\boldsymbol{V}_h:=&\;\{\boldsymbol{v}\in\boldsymbol{L}^2(\Oh): \;\boldsymbol{v}|_T\in H(\text{div}; T), \; T\in\Oh\}, 
\\
M_h:=&\;\{\mu\in{L}^2(\Eh): \; \mu|_{e}\in P^0(e), \text{  for } e\in\Eh \;  \text{  and } \; \mu|_{\partial \Omega} = 0\},\\
M_H:=&\;\{\mu\in{L}^2(\mathcal{E}_H):\; \mu|_{E}\in P^0(E), \text{  for }  E \in\EH \;\;  \text{  and } \mu|_{\partial \Omega} = 0\},\\
M^0_h:=&\;\{\mu\in{L}^2(\Eh): \; \mu|_{e}\in P^0(e), \text{  for } e\in\Eh^0 \;  \text{  and } \; \mu|_{\partial T} = 0\},\\
M_H^f:=&\;\{\mu\in M_h,\mu \notin M_h^0\},\\
M_{h,H}:=&\;M^0_h \oplus M_H,
\end{alignat*}
where $P^0(S)$ is the space of constant functions defined on the set $S$.
Note that $M_H$ is a subspace of $M_H^f$.
Then, the fine scale mortar finite element method reads as: find $\boldsymbol{q_h}\in \boldsymbol{V}_h$, $u_h\in W_h$, and $\lambda_{h,H}\in M_{h,H}$ 
such that,
\begin{subequations}\label{Eq:weak}
\begin{alignat}{2} 
\label{fine_weak-1}
(\kappa^{-1} \boldsymbol{q_h}, \boldsymbol{v})_{\mathcal{T}_h}-(u_h, \nabla\cdot\boldsymbol{v})_{\mathcal{T}_h}+\left\langle\lambda_{h,H}, \boldsymbol{v}\cdot \boldsymbol{n}\right\rangle_{\partial\mathcal{T}_h}&= 0 \qquad && \boldsymbol{v}\in \boldsymbol{V}_h\\
\label{fine_weak-2}
\left\langle\nabla\cdot \boldsymbol{q_h}, w\right\rangle_{\partial\mathcal{T}_h}&= (f,w)_{\mathcal{T}_h}\qquad&& w\in W_h\\
\label{fine_weak-3}
\left\langle \boldsymbol{q_h}\cdot \boldsymbol{n}, \mu\right\rangle_{\partial\mathcal{T}_h} &= 0\qquad && \mu\in M_{h,H}
\end{alignat}
\end{subequations}
where $\boldsymbol{n}$ is the outer normal unit vector. The inner product $(\eta, \xi)_{\mathcal{T}_h}= \sum_{T\in \mathcal{T}_h}(\eta, \xi)_T$, where $(\eta, \xi)_T$ denotes the $L^2$ inner product of $\eta$ and $\xi$ over the domain $T$, $\langle\eta, \xi\rangle_{\partial \mathcal{T}_h}= \sum_{T\in \mathcal{T}_h}\langle\eta, \xi\rangle_{\partial T}$, where $\langle\eta, \xi\rangle_{\partial T}$ denotes the $L^2$ inner product of $\eta$ and $\xi$ over the boundary $\partial T$ of the domain $T$. 

\subsection{Domain decomposition formulation}
The mortar multiscale method and the preconditioner developed in this paper
are based on a domain decomposition formulation of (\ref{fine_weak-1})-(\ref{fine_weak-3}).
We will present the mortar multiscale method in this section, and derive the preconditioner in Section \ref{two-level}.
The main feature of the mortar mixed finite element method is that it could be implemented by just solving a global system on the coarse mesh together with the solutions of some 
local problems. To achieve this, we split (\ref{fine_weak-3}) into two equations by testing separately with $\mu\in M_h^0$ and $\mu\in M_H$ so that
\begin{equation}
\label{eq:split_system}
\left\langle \boldsymbol{q_h}\cdot \boldsymbol{n}, \mu\right\rangle_{\partial\mathcal{T}_h} = 0\qquad  \mu\in M_{h}^0,
\end{equation}
and
\begin{equation}
\left\langle \boldsymbol{q_h}\cdot \boldsymbol{n}, \mu\right\rangle_{\partial K_i} = 0\qquad  \mu\in M_{H}.
\end{equation}
We can implement the solution of Eq.(\ref{eq:split_system}) independently on each subdomain $K_i$. More specifically, 
for a particular subdomain $K_i$, let $\lambda_{h,H}=\xi_H\in M_H$. Then we can find the solution $(u_h, \boldsymbol{q_h},\lambda_{h,H})|_{K_i}$
by restricting Eq.(\ref{Eq:weak}) to $K_i$:

\begin{subequations} \label{eq:split_system_2}
\begin{alignat}{2} 
(\kappa^{-1}\boldsymbol{q_h}, \boldsymbol{v})_{\mathcal{T}_h(K_i)}-(u_h, \nabla\cdot\boldsymbol{v})_{\mathcal{T}_h(K_i)}+\left\langle\lambda_{h,H}, \boldsymbol{v}\cdot \boldsymbol{n}\right\rangle_{\partial\mathcal{T}_h(K_i)}&= 0\\
\left\langle\nabla\cdot \boldsymbol{q_h}, w\right\rangle_{\partial\mathcal{T}_h(K_i)}&= (f,w)_{\mathcal{T}_h(K_i)}\\
\label{eq:split_equation_2}
\left\langle \boldsymbol{q_h}\cdot \boldsymbol{n}, \mu\right\rangle_{\partial\mathcal{T}_h(K_i)} &= 0\\
\lambda_{h,H} &= \xi_H
\end{alignat}
\end{subequations}
for all $(w, \boldsymbol{v}, \mu) \in W_h|_{\Omega_i} \times \boldsymbol{V}_h|_{K_i} \times M^0_h|_{\mathcal{E}^0_h(K_i)}$. Using the superposition principle
the solution of Eq. (\ref{eq:split_equation_2}) can be further split into two parts, namely,
\begin{equation}
(\boldsymbol{q}_h, u_h, \lambda_{h,H}) =(\boldsymbol{q}_h(f), u_h(f), 
\lambda_{h,H}(f)) 
+(\boldsymbol{q}_h(\xi_H), u_h(\xi_H), \lambda_{h,H}(\xi_H)),
\end{equation}
where $ (\boldsymbol{q}_h(f), u_h(f), \lambda_{h,H}(f)) $ satisfies
\begin{subequations}\label{eq:split_equation_a}
\begin{alignat}{2} 
(\kappa^{-1}\boldsymbol{q_h}(f), \boldsymbol{v})_{\mathcal{T}_h(K_i)}-(u_h(f), \nabla\cdot\boldsymbol{v})_{\mathcal{T}_h(K_i)}+\left\langle\lambda_{h,H}(f), \boldsymbol{v}\cdot \boldsymbol{n}\right\rangle_{\partial\mathcal{T}_h(K_i)}&= 0\\
\left\langle\nabla\cdot \boldsymbol{q_h}(f), w\right\rangle_{\partial\mathcal{T}_h(K_i)}&= (f,w)_{\mathcal{T}_h(K_i)}\\
\left\langle \boldsymbol{q_h}(f)\cdot \boldsymbol{n}, \mu\right\rangle_{\partial\mathcal{T}_h(K_i)} &= 0\\
\lambda_{h,H}(f)&= 0
\end{alignat}
\end{subequations}
for all $(w, \boldsymbol{v}, \mu) \in W_h|_{K_i} \times \boldsymbol{V}_h|_{K_i} \times M^0_h|_{\mathcal{E}^0_h(K_i) }$
and $ (\boldsymbol{q}_h(\xi_H), u_h(\xi_H), \lambda_{h,H}(\xi_H))$ satisfies

\begin{subequations}\label{eq:split_equation_b}
\begin{alignat}{2} 
(\kappa^{-1} \boldsymbol{q_h}(\xi_H), \boldsymbol{v})_{\mathcal{T}_h(K_i)}-(u_h(\xi_H), \nabla\cdot\boldsymbol{v})_{\mathcal{T}_h(K_i)}+\left\langle\lambda_{h,H}(\xi_H), \boldsymbol{v}\cdot \boldsymbol{n}\right\rangle_{\partial\mathcal{T}_h(K_i)}&= 0\\
\left\langle\nabla\cdot \boldsymbol{q_h}(\xi_H), w\right\rangle_{\partial\mathcal{T}_h(K_i)}&= 0\\
\left\langle \boldsymbol{q_h}(\xi_H)\cdot \boldsymbol{n}, \mu\right\rangle_{\partial\mathcal{T}_h(K_i)} &= 0\\
\lambda_{h,H}(\xi_H)&= \xi_H
\end{alignat}
\end{subequations}
for all $(w, \boldsymbol{v}, \mu) \in W_h|_{K_i} \times \boldsymbol{V}_h|_{K_i} \times M^0_h|_{\mathcal{E}^0_h(K_i)}$.

Then Eq. (\ref{eq:split_system_2}) can be reduced to finding $\xi\in M_H$ such that
\begin{equation}
\label{coarse-system}
a_H(\xi_H,\mu)=g_H(\mu)\quad\quad \text{for all } \mu\in M_H,
\end{equation}
where the bilinear form $a_H(\xi_H, \mu): M_H\times M_H\to R$
and the linear form $ g_H(\mu): M_H \to R$ are defined as
\begin{equation}\label{form-a}
a_H(\xi_H, \mu):=\sum_{i=1}^{N}\left\langle\lambda_{h,H}(\xi_H), \mu\right\rangle_{\partial K_i} \quad \text{and} \quad  g_H(\mu):=-\sum_{i=1}^{N}\left\langle\lambda_{h,H}(f), \mu\right\rangle_{\partial K_i} . 
\end{equation} 
The interface bilinear form $a_H(\cdot,\cdot)$ is symmetric and positive semi-definite on $M_H$ and this system can be solved by preconditioned conjugate gradient method. See \cite{Arbogast_PWY_07,hdgunified,balancing2003} and reference therein for more details. We will construct a preconditoner for (\ref{coarse-system}) in Section \ref{two-level}.
On the other hand, the mortar space $M_H$ in (\ref{coarse-system}) is not accurate enough for problems with high contrast heterogeneous media.
Thus, in the next section, we will introduce a multiscale mortar space for (\ref{coarse-system}).
In particular, we will replace the space $M_H$ by our multiscale mortar space.

\section{Multiscale mortar space}\label{ms-basis}
In this section, we present the construction of our multiscale mortar space. 
The space can be used as an approximation space for a multiscale mortar method in (\ref{coarse-system}),
and as a coarse solver for a class of preconditioners, which will be presented in Section \ref{two-level}.
Our multiscale mortar space consists of a set of multiscale basis functions,
which are defined only on the coarse skeleton ${\EH}$. 
To construct these basis functions, we will first define a set of snapshot functions for each coarse edge. 
The snapshots represent various modes of the solutions, and are typically large. 
To find the snapshots, we will solve some local problems in a small subdomain covering an edge,
and then restrict the solutions to the edge. 
We will next define some PODs and use them to extract dominant modes within the snapshot space.
These dominant modes form the multiscale basis functions. 
Notice that these basis functions capture some information of the permeability field within neighboring coarse blocks of an edge. 


\textbf{Construction of local snapshot space} $V^E_{\text{snap}}(\omega_i)$.

Let $E_i \in {\cal E}_H^0$ be an interior edge and let $\omega_i$ be the corresponding coarse neighborhood covering $E_i$ (see Figure~\ref{grids}). 
We will first construct a set of
local snapshots $\{\psi_j^{E_i}\}_{j=1}^{N_i}$  defined on $\omega_i$ and $$V^{E_i}_{\text{snap}}(\omega_i)=\text{span}\{\psi_j^{E_i}\}_{j=1}^{N_i}.$$
The snapshots can be given explicitly or computed via solutions of local boundary value or local spectral problems in $\omega_i$.
We will in this paper use the following approach. For each coarse edge $E_i$, we define
$$W_i(\partial \omega_i)=\{ w_j^i \; | \; w_j^i= 1 \quad \mbox{on}\quad e_j; \quad w_j^i= 0 \quad \mbox{on} \quad \partial \omega_i-e_j, 1\leq j\leq N_{\omega_i}\},$$ where
$N_{\omega_i}$ is the number of fine eges on $\partial \omega_i$, and $e_j$ is the $j$-th fine edge on $\partial\omega_i$.
To construct the snapshots, we solve the following problem 
\begin{eqnarray}
-\nabla\cdot(\kappa\nabla u_j^{i})&=&0,\quad\quad \text{in} \quad \omega_i,\\
u_j^{i}&=&w_j^i,\quad \quad \text{on}\quad \partial\omega_i.
\end{eqnarray}
Using the solutions of the above problem, we obtain
the space $\{ \psi_j^E=u_j^{i}|_{E_i}, j=1,\cdots N_{\omega_i}\}. $ 
Finally, we obtain the snapshot space $V^E_{\text{snap}}(\omega_i)=\text{span} \{\psi_j^E\}.$

\textbf{Construction of multiscale mortar space} $V^E_{\text{off}}$. 

To construct our multiscale mortar space, we will apply a space reduction technique to the snapshot space $V^{E_i}_{\text{snap}}(\omega_i)$
to obtain a smaller dimensional space. 
In particular, we will perform POD to $V^E_{\text{snap}}(\omega_i)$ and then select the first $l_i$ dominant modes $\Psi_j^{i}$.
In this way, we obtain the offline space corresponding to the coarse face $E_i$, which is
\begin{equation}
\label{off}
V^E(E_i)=\mbox{span}\{\Psi_j^{i}, 1\leq j\leq l_i \}.
\end{equation}
To simplify the notations, we use the single-index notation:
$$V^E_{\textrm{off}}=\textrm{span}\{ \Psi_i^{\textrm{off}}:  1\leq i \leq M_{\textrm{off}}\},$$
where $M_{\textrm{off}}=\sum_{i=1}^{N_e}J_i.$
We define $$ R_{\textrm{off}}^E=[\psi_1^{\textrm{off}},\ldots,\psi_{M_{\textrm{off}}}^{\textrm{off}}],$$ which maps from the offline space to the fine space, and
$\psi_i^{\textrm{off}}$ is a vector containing the coefficients in the expansion of $\Psi_i^{\textrm{off}}$ in the fine-grid basis functions.
We use this multiscale  offline basis to enrich the constant mortar space $M_H$,
and use the resulting space to solve (\ref{coarse-system}).

\textbf{Oversampling and randomized snapshot}.

One can use an oversampling technique to improve the accuracy of the method. The main idea
of the oversampling approach is to use  larger domains $\omega_i^+$ instead of $\omega_i$ to compute snapshot, see Figure~\ref{oversampgrid} for examples of oversampling domain of a coarse edge.  {A typical choice of $V^E_{\text{snap}}(\omega_i)$ is local flow solutions, which are constructed by solving local problems with all possible boundary conditions. This can be expensive and for this reason, one can use randomized boundary conditions (see \cite{calo2014randomized}) and construct only a few more snapshots than the number of desired local basis functions.

\section{Two-level domain decomposition preconditioners}\label{two-level}
Previous sections are devoted to the design a good multiscale mortar space to achieve a desired accuracy. 
In this section, we will apply our multiscale mortar space for an
iterative method applied to solve
the fine scale system. We will describe a two-level iterative method to solve the fine scale problem with the previously introduced 
multiscale basis as the coarse space, see \cite {arbogast2015two,xiao2013multiscale} for the case of polynomial and homogenized multiscale basis.
The two-level preconditioner $B^{-1}$ 
includes two part, a local fine scale preconditioner $B_{loc}^{-1}$  to smooth out fine-scale error and a global coarse-scale preconditioner $B_0^{-1}$ for exchanging global information. 

\subsection{Global coarse preconditioners}
We denote the coarse basis functions as $\{\Phi_i\}^{N_c}_
{i=1}$,
where $N_c$ the total number of coarse basis functions, which is usually larger than
the total number of coarse edges due to the use of enriched basis. Then we can
define the coarse space as
\begin{equation*}
V_0=\text{span}\{\Phi_i\}^{N_c}_{i=1}
\end{equation*} 
and the coarse matrix $A_0=R_0AR_0^T$,
where $R_0^T= [\Phi_1,\Phi_2,...,\Phi_{N_c}]$, and $A$ is the matrix corresponding to bilinear form
$a_H(\xi_i, \xi_j): M_H^f\times M_H^f\to R$ defined in  (\ref{form-a}). Then, the coarse preconditioner is defined 
as 
\begin{equation*}
 B_0^{-1}=R_0A_0^{-1}R_0^T.
 \end{equation*}
 The selection of coarse space $V_0$ is quite important to the performance of the preconditioner, 
we will use the multiscale basis (\ref{off}) constructed 
in the previous section to form the mortar space $V_0$.

\subsection{Local preconditioners}
The local preconditioner is $B_{\text{loc}}^{-1}$ defined neighborhood wise. More specifically, 
Let $\mathcal R_i: M_H\to M_H|_{E_i}$ be the restriction operator from $\EH$ to $E_i$
and let $R_i$ be the corresponding matrix representation. For each coarse edge $E_i$, we will consider 
 a domain $\omega_i^+\supset E_i$ (see Figure~\ref{oversampgrid} for the illustration
of $\omega_i^{+}$) to apply the local preconditioner. Similarly, we define the restriction operator from $\EH$ to $E_i^+=\EH\cap\omega_i^{+} $ as $\mathcal P_i: M_H\to M_H|_{E_i^+}$
and corresponding matrix $P_i$. Note that $\omega_i^+$ can be the same as $\omega_i$.
Then we  define the local preconditioner as
\begin{equation*}
B_{\text{loc}}^{-1}=\sum\limits_{i} R_{i}^TA_{i}^{-1}P_{i}, 
\end{equation*}
where $A_i=P_iAP_i^{T}$.
The application of $A_{i}^{-1}$ is equivalent to solve a homogeneous Dirichlet boundary 
condition problem on $\omega_i^+$ with local residual as source.  

If $\omega_i^+=\omega_i$, then we have $R_{i}=P_i$, which implies that $B_{\text{loc}}^{-1}$ is symmetric
and we can use PCG as outside accelerator. In this case, the computational
cost (although it is offline) of applying local preconditioners may be expensive especially
in 3D case. 
To reduce the computational cost, we can
consider the case $\omega_i^+\neq\omega_i$, which is actually
the restrictive local preconditioner \cite{cai1999restricted}. 
 This will not only reduce the computation of applying local preconditioner, but 
also decrease the number of iterations since it includes the distant information.
In this case $B_{\text{loc}}^{-1}$ is no longer symmetric, and we can choose algorithm such as
GMRES  as the outer accelerator.

\begin{remark}
	Restrictive local preconditioner is quite similar with the idea of oversampling, both
	utilize a larger domain than standard domain to perform computation and then take restriction. Both method shows better performance than standard method.
\end{remark}

\subsection{Two-level preconditioners}
 We combine the local preconditioner and coarse preconditioner 
in two ways to form the two-level preconditioners. The first approach is the {\it additive} preconditioner
\begin{equation}
B_{\text{add}}^{-1}= B_0^{-1}+ B_{\text{loc}}^{-1}.
\end{equation}
 The second is the {\it hybrid} preconditioner
\begin{equation}
B_{\text{hyb}}^{-1}= B_0^{-1}+ (I-B_0)B_{\text{loc}}^{-1}(I-B^T_0).
\end{equation}
where $P_0=B_0A$ is the Schwarz projection operator, see \cite{toselli2005domain} for details.
For more details about the two-level preconditioners we adopted here, we refer \cite {arbogast2015two,xiao2013multiscale} and reference therein.
We remark again that the main ingredient in the above preconditioners
is the use of our mortar multiscale space constructed in Section \ref{ms-basis}.

\section{Numerical examples}\label{numerical-results}
In this section, we present some representative examples to show the performance of our method. We consider two models with permeability $\kappa$ depicted in Figure~\ref{model}. We note that $\kappa=1$ in the blue region and $\kappa=\eta (>>1)$ in the red region.  We will consider two cases: $\eta=10^4, \eta=10^6$ in the following examples. As it is shown, these two models contains high contrast, short and long channels, and isolated inclusions.
We will first demonstrate the performance of the multiscale solver by showing the error of multiscale solution against  the fine scale (reference) solution. Next we report the results of  two-level additive Schwarz domain decomposition method with the coarse space formed by using our multiscale mortar space. 
We will consider different snapshot space computed on different domains, and also consider various preconditioners.

 We divide the domain $\Omega=(0,1)^2$ into $N\times N$ square coarse elements. In each coarse element, we generate a uniform $n\times n$ fine scale square elements. Therefore, the domain was divided into $N_f\times N_f$  fine elements, where $N_f =N\times n$. The number of degrees of freedom for the fine solver is $5\times N_f^2 + 2\times N_f \times (N_f-1)$, while the number of degree of freedom for the multiscale  solver is $N_b\times 2\times N \times (N-1)$, where $N_b$ is the number of  multiscale basis on each coarse edge.

Constant sources and homogeneous Dirichlet boundary condition are considered.
 We use $\left(\begin{array}{cc}
d_{11} & d_{12} \\
d_{21}  & d_{22} \\
 \end{array}
 \right)$ to define local computational domain of the snapshot and the local preconditioner. See Figure~\ref{oversampgrid} for the illustration of 
 $d_{ij}$.
In total, 4 ways to generate the multiscale space are considered:\\
Domain 1: No oversampling: $\left(\begin{array}{cc}
$n$ & 0 \\
0 & $n$\\
\end{array}
\right)\\$
Domain 2: oversampling case a: $\left(\begin{array}{cc}
$n$ & 1 \\
1& $n$\\
\end{array}
\right)\\$\\
Domain 3: oversampling case b: $\left(\begin{array}{cc}
$ [n/2] $ & 1 \\
1 & $ [n/2] $\\
\end{array}
\right)\\$\\
Domain 4: oversampling case c: $\left(\begin{array}{cc}
2 & 1 \\
1 & 2\\
\end{array}
\right)\\$

\begin{figure}[ht]
  \centering
  \subfigure[$\kappa_{1}$]{
    \includegraphics[width=3.0in]{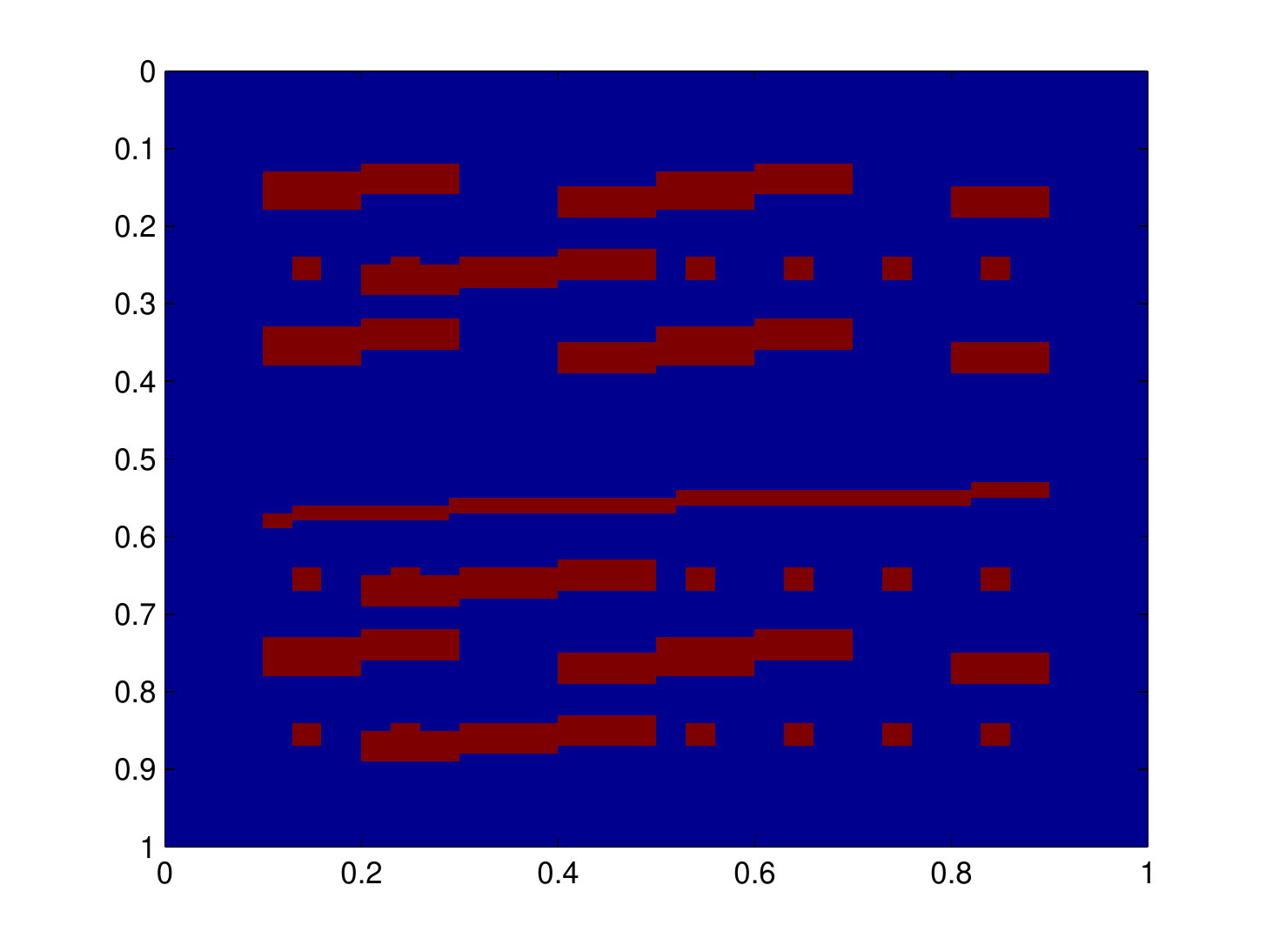}}
   \subfigure[$\kappa_{2}$]{
     \includegraphics[width=3.0in]{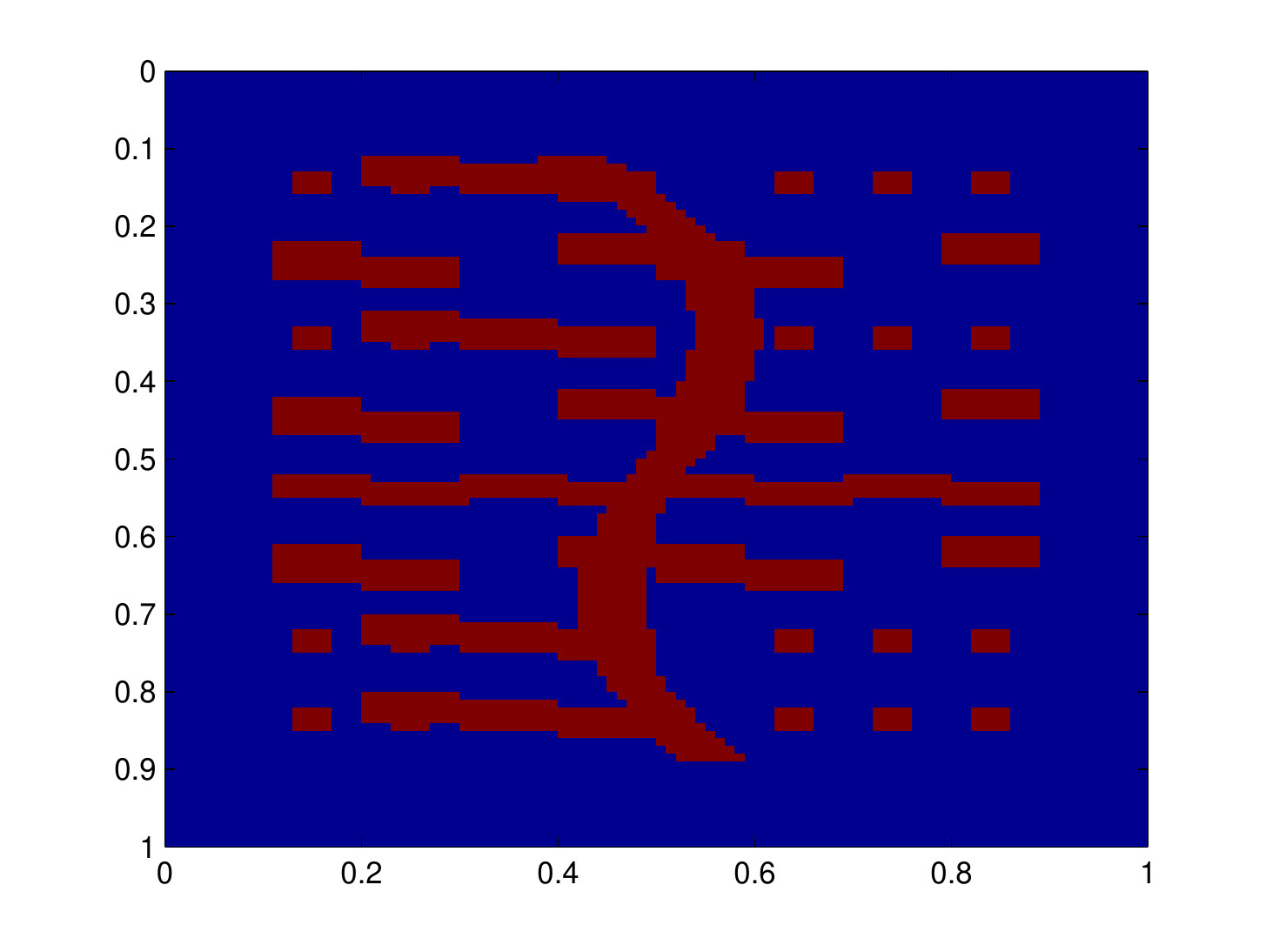}}
  \caption{Permeability fields.}
  \label{model} 
\end{figure}

\begin{figure}[ht]
\centering
\resizebox{0.7\textwidth}{!}{
\begin{tikzpicture}[scale=0.8]
\filldraw[fill=orange, draw=black] (1,0) rectangle (5,5);
\draw[step=1,black, thick] (1,0) grid (5,5);
\draw[step=1,black, thick] (0,1) grid (6,4);
\draw[ultra thick, green](3, 1) -- (3,4);
\node at (3.5,2.5)  { \huge $E_{i}$} ;
\node at (12.5,3)  { \huge $\omega_{i}^+$} ;
\draw[ thick, black](1, 5.1) -- (1,5.4);
\draw[ thick, black](1, 5.2) -- (1.6,5.2);
\node at (2,5.3)  { $d_{11}$};
\draw[ thick, black](2.4, 5.2) -- (3,5.2);
\draw[ thick, black](3, 5.1) -- (3,5.4);
\draw[ thick, black](3, 5.2) -- (3.6,5.2);
\node at (4,5.3)  { $d_{11}$};
\draw[ thick, black](4.4, 5.2) -- (5,5.2);
\draw[ thick, black](5, 5.1) -- (5,5.4);

\draw[ thick, black](6.1, 5) -- (6.4,5);
\draw[ thick, black](6.3, 5) -- (6.3,4.7);
\node at (6.3,4.6)  { $d_{12}$};
\draw[ thick, black](6.1, 4) -- (6.4,4);
\draw[ thick, black](6.3, 4.4) -- (6.3,4);
\draw[ thick, black](6.3, 4) -- (6.3,2.6);
\node at (6.3,2.5)  { $H$};
\draw[ thick, black](6.3, 2.3) -- (6.3,1);
\draw[ thick, black](6.1, 1) -- (6.4,1);
\draw[ thick, black](6.3, 1) -- (6.3,0.7);
\node at (6.3,0.6)  { $d_{12}$};
\draw[ thick, black](6.3, 0.4) -- (6.3,0);
\draw[ thick, black](6.1, 0) -- (6.4,0);

\filldraw[fill=orange, draw=black] (13,1) rectangle (18,5);
\draw[step=1,black, thick] (13,1) grid (18,5);
\draw[step=1,black, thick] (14,0) grid (17,6);

\draw[ultra thick, green](14, 3) -- (17,3);
\node at (15.5,2.5)  { \huge $E_{i}$} ;
\node at (3,-0.5)  { \huge $\omega_{i}^+$} ;
\draw[ thick, black](14,-0.1) -- (14,-0.4);
\draw[ thick, black](13,-0.1) -- (13,-0.4);
\draw[ thick, black](17,-0.1) -- (17,-0.4);
\draw[ thick, black](18,-0.1) -- (18,-0.4);
\draw[ thick, black](13, -0.2) -- (13.2,-0.2);
\node at (13.5,-0.2)  { $d_{21}$};
\draw[ thick, black](13.7, -0.2) -- (14,-0.2);
\draw[ thick, black](14, -0.2) -- (15.2,-0.2);
\draw[ thick, black](15.7, -0.2) -- (17,-0.2);
\node at (15.5,-0.2)  { $H$};
\draw[ thick, black](17, -0.2) -- (17.2,-0.2);
\draw[ thick, black](17.7, -0.2) -- (18,-0.2);
\node at (17.5,-0.2)  { $d_{21}$};

\draw[ thick, black](18.1,1) -- (18.4,1);
\draw[ thick, black](18.1,3) -- (18.4,3);
\draw[ thick, black](18.1,5) -- (18.4,5);
\draw[ thick, black](18.3,5) -- (18.3,4.1);
\draw[ thick, black](18.3,3.7) -- (18.3,3);
\node at (18.3,3.9)  { $d_{22}$};
\draw[ thick, black](18.3,3) -- (18.3,2.1);
\draw[ thick, black](18.3,1.7) -- (18.3,1);
\node at (18.3,1.9)  { $d_{22}$};
\end{tikzpicture}
}
\caption{ Illustration of a oversampled neighborhood.}
\label{oversampgrid}
\end{figure}
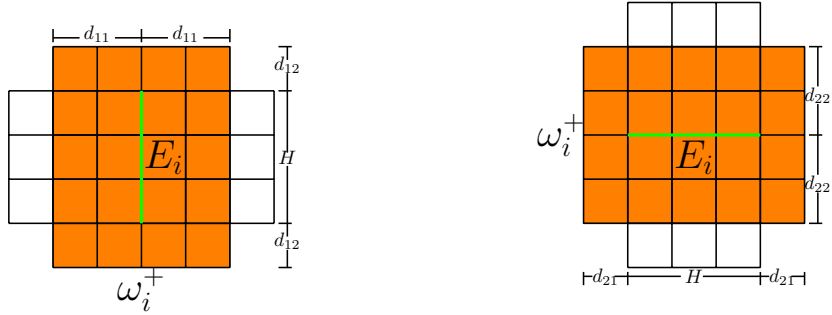

\subsection{Coarse grid multiscale solution}
In this subsection, we study the error decay of multiscale solution by adding more multiscale basis. We will consider using 4 types of multiscale basis as well as polynomial basis for comparison. We will
demonstrate the influence of snapshot, mesh size and contrast of the permeability on the multiscale solution.  We define the following error to quantify the accuracy of coarse 
grid multiscale solution.
\begin{equation*}
e_u:=\frac{\|u_{ms}-u_f\|_{L^2,\Omega}}{\|u_f\|_{L^2,\Omega}}, \quad  e_{\boldsymbol{q}}:=\frac{\|\boldsymbol{q}_{ms}-\boldsymbol{q}_f\|_{\kappa,\Omega}}{\|\boldsymbol{q}_f\|_{\kappa,\Omega}}
\end{equation*}
where $\|\boldsymbol{q}\|_{\kappa,\Omega}^2=\int_{\Omega}\kappa^{-1}\boldsymbol{q}^2dx$.

The 4 types of multiscale basis based on  the way of generating the snapshot are:\\
Case 1: full snapshot on domain 1\\
Case 2: full snapshot  on domain 2\\
Case 3: oversampling full snapshot on domain 3\\
Case 4: $(n+2)$ randomized snapshot on domain 2.

Table \ref{model1_elliptic_5} and Table \ref{model1_elliptic_10} present the numerical results  for model 1 with  mesh setting $N=5, n=20$ and $N=10, n=10$ respectively. The first column is the number of multiscale basis for each coarse edge. The rest of the columns present both errors $e_u, e_{\boldsymbol{q}}$ for using polynomial basis, and 4 types of multiscale basis as described before. As it is shown, by adding basis, both types of errors decrease for all types of basis considered. Moreover, the error decay of all multiscale basis cases are faster than the polynomial basis case.
 For example, in Table \ref{model1_elliptic_5}, the  relative weighted velocity error $e_{\boldsymbol{q}}$ decreases from 78.6\% to 6.2\% in case 2 when $N=5$; however, the corresponding error decreases from 78.2\% to 26.7\% for  polynomial basis. We remark that 6\% error is tolerable in industry flow simulation.  By using 5  basis on each coarse edge, the error $e_{u}$ for the polynomial case is 10.3\%, while this error for the multiscale basis cases are below ten percent. Note that, by using 5  basis on each coarse edge, the number of degree of freedom of multiscale solver is about 0.3\% of that of the fine solver.
 By comparing case 1 and case 2(oversampling case), 
we can observe obvious improvement by applying oversampling although only 1 fine scale element is added to the snapshot domain. We remark that adding more fine scale elements to the domain of computing snapshot will
further improve the solution. From the results of case 2 and case 4, we can see that the performance of randomized snapshot is also comparable with full snapshot, and it is better than case 3 and case 1. In Figure~\ref{model1_pressure}, for model 1, the reference solution and the corresponding multiscale solutions
with different number of basis on each coarse edge are presented.  The upper left is the reference solution. The upper right is  the multiscale solution with 1 basis on each coarse edge. This solution has obvious discontinuity across the coarse edges. The bottom left shows the multiscale solution with 3 basis on each coarse edge, which is closer to the reference solution. However, we can see slight difference between this solution and reference solution. The bottom  right is the multiscale solution with 5 basis on each coarse edge,
this solution is almost identical with the reference solution. This figure shows that the additional multiscale basis functions are important to capture all the features of the solution.

By comparing the errors in Table \ref{model1_elliptic_5} and Table \ref{model1_elliptic_10}, we notice that smaller coarse grid size can improve the accuracy of the solution. For example, the second column in Table \ref{model1_elliptic_5}, $e_u$ for polynomial basis drops from 63.0\% to 10.3\%, while the same error in Table \ref{model1_elliptic_10} drops from 60.9\% to 5.0\%; the third column in Table \ref{model1_elliptic_5}, $ e_{\boldsymbol{q}}$ for Case 1 drops from 78.6\% to 8.8\%, while the same error in Table \ref{model1_elliptic_10} drops from 77.8\% to 2.7\%.

Table \ref{model2_elliptic_5} and Table \ref{model2_elliptic_10} show corresponding results for model 2. We observe similar results as model 1.  Both types of errors decrease for all types of basis considered by adding more basis. The error decay of all multiscale basis cases are faster than the polynomial basis case. For example, if $N=10$  by adding to 5 basis on each coarse edge, the error $e_{u}$ for polynomial case decreases to 1.9\%; while the errors for multiscale cases decrease to 0.6\%, 0.04\%, 1.1\%, 0.08\% respectively. 
 By comparing case 1 and case 2, we can again observe obvious improvement by applying oversampling although only 1 fine scale element is added to the snapshot domain. From the results of case 2 and case 4, we can see that the performance of randomized snapshot is also comparable with full snapshot, and it is better than case 3 and case 1.

In Figure~\ref{model2_pressure}, for model 2, the reference solution and the corresponding multiscale solutions
with different numbers of basis on each coarse edge are presented.  The upper left is the reference solution. The upper right is  the multiscale solution with 1 basis on each coarse edge. Apparently the solution is discontinuous across the coarse edges. The bottom left is the multiscale solution with 3 basis on each coarse edge, which shows better agreement with the reference solution. However, there are still some detailed features missing in the multiscale solution . The bottom  right is the multiscale solution with 5 basis on each coarse edge, which has good agreement  with the reference solution. This figure shows that the additional multiscale basis functions are important to capture all the features of the solution. These results for model 2 demonstrate that our multiscale solver can handle permeability field with long channel effects.

We also test the robustness of our method by varying the order of high contrast. The results are presented in Figure~\ref{contrast_model1} and Figure~\ref{contrast_model2}. In Figure~\ref{contrast_model1}, the errors $e_u$(left), $e_{\boldsymbol{q}}$ (right) for model 1 with high contrast order $\eta= 10^4, 10^6$  with coarse $N=5, N=10 $ are plotted. The black dashed line is for the case $10^4, N=5$, the blue dashed line is for the case $10^6, N=5$. The two lines are identical. The green dashed line is for the case $10^4, N=10$, the red dashed line is for the case $10^6, N=10$. The two lines are also identical. 
 
Figure~\ref{contrast_model2} displays the corresponding results for model 2. Though the lines do not overlap, the difference is very small. We note that the snapshot we used for both cases comes from case 2. From these examples, we see that our method produces robust results independent of the order of high contrast. Robustness will further be confirmed by our preconditioner
results in following parts. 

\begin{table}[H]
\centering
\begin{tabular}{|c|c|c|c|c|c|c|c|c|c|c|}
\hline 
\multirow{2}{*}{Nb} & \multicolumn{2}{c|}{\specialcell{Polynomial}}& \multicolumn{2}{c|}{\specialcell{Case 1}} & \multicolumn{2}{c|}{ \specialcell{Case 2}} & \multicolumn{2}{c|}{\specialcell{Case 3}}& \multicolumn{2}{c|}{\specialcell{Case 4}}\tabularnewline
\cline{2-11} 
 & $e_u$ &$e_{\boldsymbol{q}}$ & $e_u$ &$e_{\boldsymbol{q}}$ & $e_u$ &$e_{\boldsymbol{q}}$ & $e_u$ & $e_{\boldsymbol{q}}$& $e_u$ & $e_{\boldsymbol{q}}$ \tabularnewline
\hline
1&0.630& 0.786& 0.630& 0.786&0.630& 0.786&0.630& 0.786&0.630& 0.786  \tabularnewline
\hline                                                                                
2 &0.432&0.647&0.438 &0.649 &0.357&0.581&0.351& 0.581&0.427& 0.640 \tabularnewline
\hline
3&0.288&0.519&0.157 & 0.374&0.139&0.354 &0.143&0.357&0.129&0.341  \tabularnewline
\hline
4&0.120&0.291& 0.060& 0.232&0.014&0.104&0.050&0.207  &0.038&0.171  \tabularnewline
\hline
5&0.103&0.267 & 0.012&0.088 &0.006& 0.062&0.028& 0.157 &0.008&0.073 \tabularnewline
\hline 
\end{tabular}
\caption{Relative error between multiscale solution and fine scale solution with different types of basis for model 1, $N=5$, $\eta=10^4$. "Nb" represent the number of basis per coarse edge.}
\label{model1_elliptic_5}
\end{table}

\begin{table}[H]
\centering
\begin{tabular}{|c|c|c|c|c|c|c|c|c|c|c|}
\hline 
\multirow{2}{*}{Nb} & \multicolumn{2}{c|}{\specialcell{Polynomial}}& \multicolumn{2}{c|}{\specialcell{Case 1}} & \multicolumn{2}{c|}{ \specialcell{Case 2}} & \multicolumn{2}{c|}{\specialcell{Case 3}}& \multicolumn{2}{c|}{\specialcell{Case 4}}\tabularnewline
\cline{2-11} 
 & $e_u$ &$e_{\boldsymbol{q}}$ & $e_u$ &$e_{\boldsymbol{q}}$ & $e_u$ &$e_{\boldsymbol{q}}$ & $e_u$ & $e_{\boldsymbol{q}}$& $e_u$ & $e_{\boldsymbol{q}}$ \tabularnewline
\hline
1&0.609 &0.778&0.609 &0.778 &0.609&0.778 &0.609 &0.778&0.609 &0.778  \tabularnewline
\hline
2 &0.133&0.322&0.140 &0.329&0.130&0.316&0.133&0.320 &0.106& 0.287 \tabularnewline
\hline
3&0.103&0.278&0.026 & 0.128&0.019&0.114&0.039&0.167&0.026& 0.135 \tabularnewline
\hline
4&0.069&0.232& 0.012&0.092 &0.010& 0.083&0.029& 0.142 &0.004&0.050  \tabularnewline
\hline
5&0.050&0.197 & 0.001&0.027 &3.9e-04& 0.013&0.002&0.035  &6.7e-04&0.020 \tabularnewline
\hline 
\end{tabular}
\caption{Relative error between multiscale solution and fine scale solution with different types of basis for model 1, $N=10$, $\eta=10^4$. "Nb" represent the number of basis per coarse edge.}
\label{model1_elliptic_10}
\end{table}

\begin{table}[H]
\centering
\begin{tabular}{|c|c|c|c|c|c|c|c|c|c|c|}
\hline 
\multirow{2}{*}{Nb} & \multicolumn{2}{c|}{\specialcell{Polynomial}}& \multicolumn{2}{c|}{\specialcell{Case 1}} & \multicolumn{2}{c|}{ \specialcell{Case 2}} & \multicolumn{2}{c|}{\specialcell{Case 3}}& \multicolumn{2}{c|}{\specialcell{Case 4}}\tabularnewline
\cline{2-11} 
 & $e_u$ &$e_{\boldsymbol{q}}$ & $e_u$ &$e_{\boldsymbol{q}}$ & $e_u$ &$e_{\boldsymbol{q}}$ & $e_u$ & $e_{\boldsymbol{q}}$& $e_u$ & $e_{\boldsymbol{q}}$ \tabularnewline
\hline
1&0.668&0.806&0.668&0.806&0.668&0.806&0.668&0.806&0.668&0.806 \tabularnewline
\hline
2 &0.508&0.699&0.549& 0.728&0.410 &0.627  &0.447&0.655&0.466&0.651 \tabularnewline
\hline
3&0.332&0.560 &0.265&0.494&0.151 &0.370 &0.151&0.370&0.142&0.358\tabularnewline
\hline
4&0.118&0.316&0.136&0.348&0.038 &0.178  &0.039&0.181&0.041&0.187 \tabularnewline
\hline
5 &0.077&0.244 &0.061&0.211&0.012 &0.089&0.009&0.075&0.026&0.153\tabularnewline
\hline 
\end{tabular}
\caption{Relative error between multiscale solution and fine scale solution with different types of basis for model 2, $N=5$, $\eta=10^4$. "Nb" represent the number of basis per coarse edge.}
\label{model2_elliptic_5}
\end{table}

\begin{table}[H]
\centering
\begin{tabular}{|c|c|c|c|c|c|c|c|c|c|c|}
\hline 
\multirow{2}{*}{Nb} & \multicolumn{2}{c|}{\specialcell{Polynomial}}& \multicolumn{2}{c|}{\specialcell{Case 1}} & \multicolumn{2}{c|}{ \specialcell{Case 2}} & \multicolumn{2}{c|}{\specialcell{Case 3}}& \multicolumn{2}{c|}{\specialcell{Case 4}}\tabularnewline
\cline{2-11} 
 & $e_u$ &$e_{\boldsymbol{q}}$ & $e_u$ &$e_{\boldsymbol{q}}$ & $e_u$ &$e_{\boldsymbol{q}}$ & $e_u$ & $e_{\boldsymbol{q}}$& $e_u$ & $e_{\boldsymbol{q}}$ \tabularnewline
\hline
1&0.593&0.769&0.593&0.769&0.593& 0.769&0.593&0.769&0.593&0.769 \tabularnewline
\hline
2 &0.055&0.205&0.118&0.304&0.090&0.266&0.090&0.267&0.074&0.245 \tabularnewline
\hline
3&0.044&0.175&0.031&0.146&0.011&0.082&0.017&0.109&0.014&0.098\tabularnewline
\hline
4&0.029&0.142&0.015&0.101&0.002&0.035 &0.007&0.065&0.004&0.049  \tabularnewline
\hline
5 &0.019&0.115&0.006&0.055&4.4e-04&0.011 &0.001&0.023&7.6e-04&0.017 \tabularnewline
\hline 
\end{tabular}
\caption{Relative error between multiscale solution and fine scale solution with different types of basis for model 2, $N=10$, $\eta=10^4$. "Nb" represent the number of basis per coarse edge.}
\label{model2_elliptic_10}
\end{table}

\begin{figure}[H]
	\centering
	\subfigure[ Fine-scale solution]{
		\includegraphics[width=3.0in]{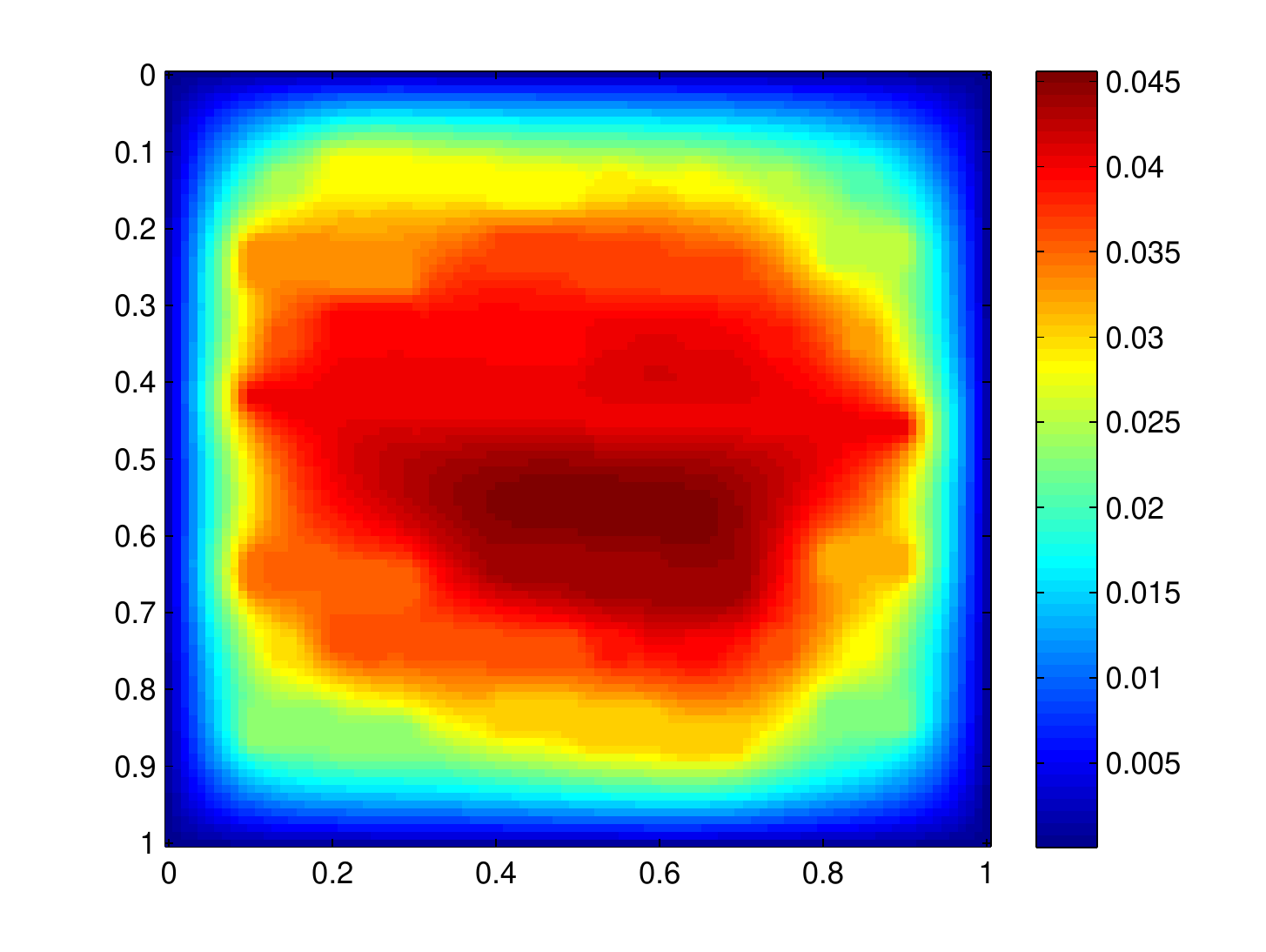}}
	\subfigure[Coarse-scale solution(1 basis)]{
		\includegraphics[width=3.0in]{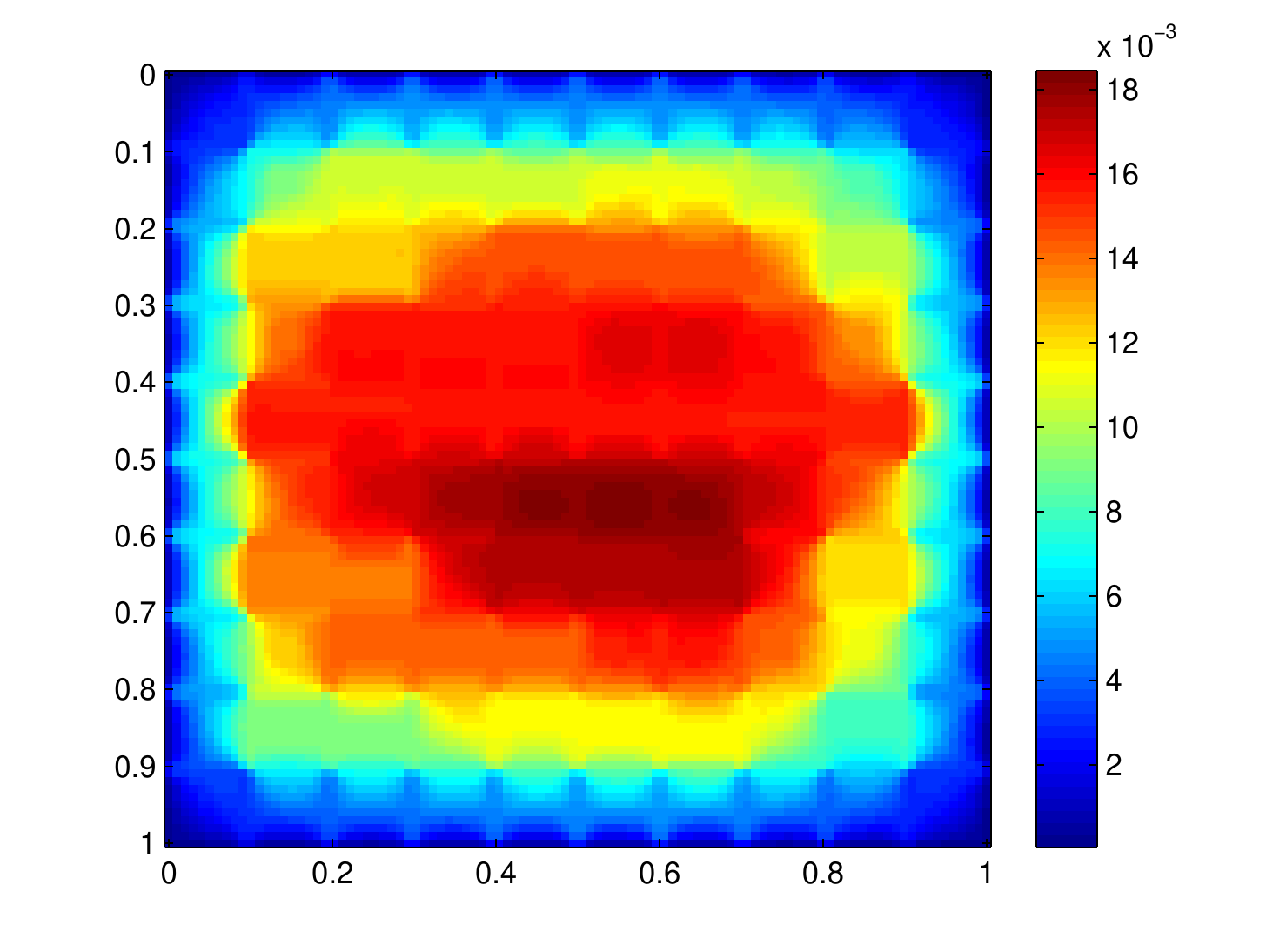}}
		\subfigure[Coarse-scale solution(3 basis)]{
			\includegraphics[width=3.0in]{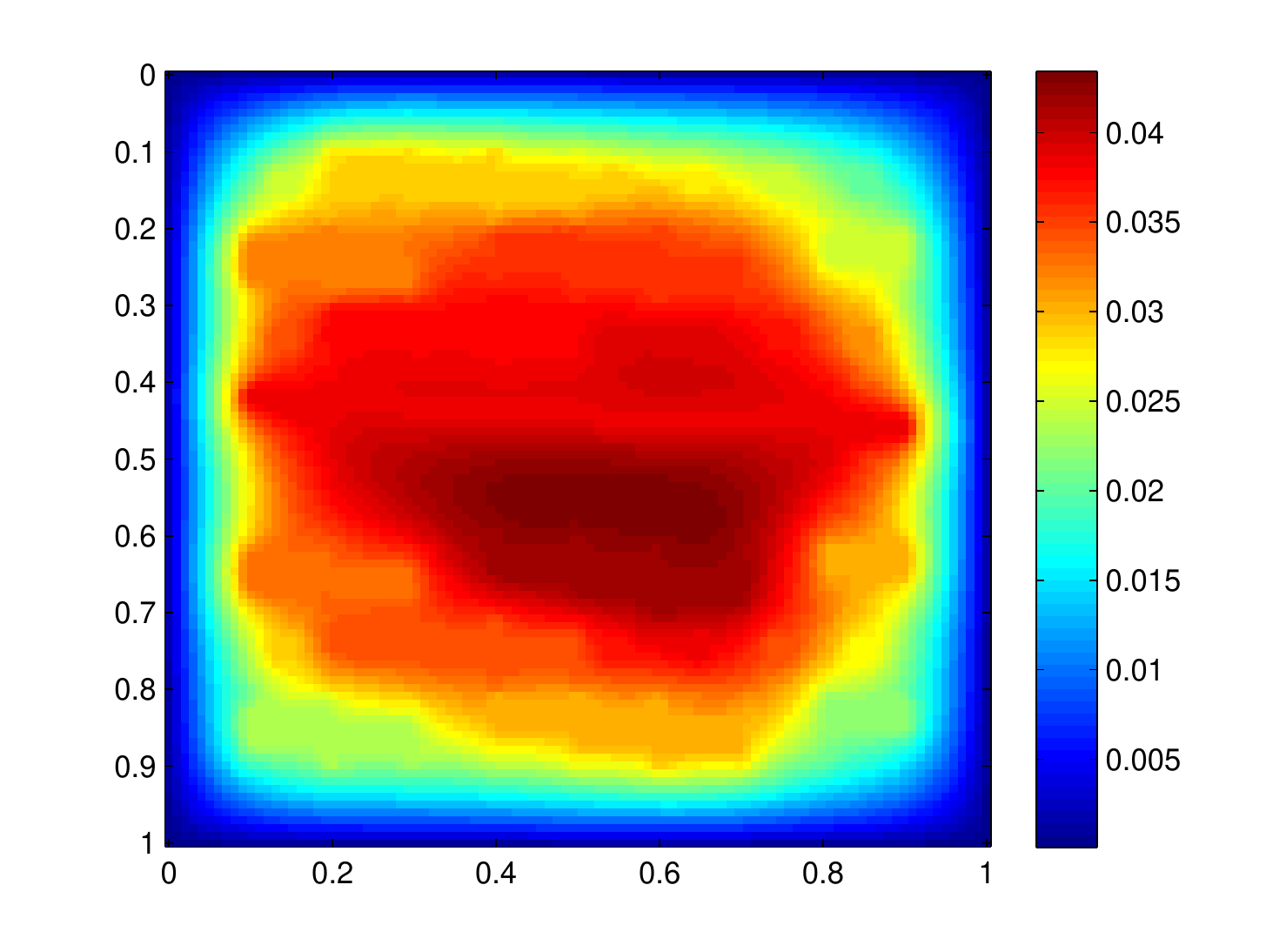}}
		\subfigure[Coarse-scale solution(5 basis)]{
			\includegraphics[width=3.0in]{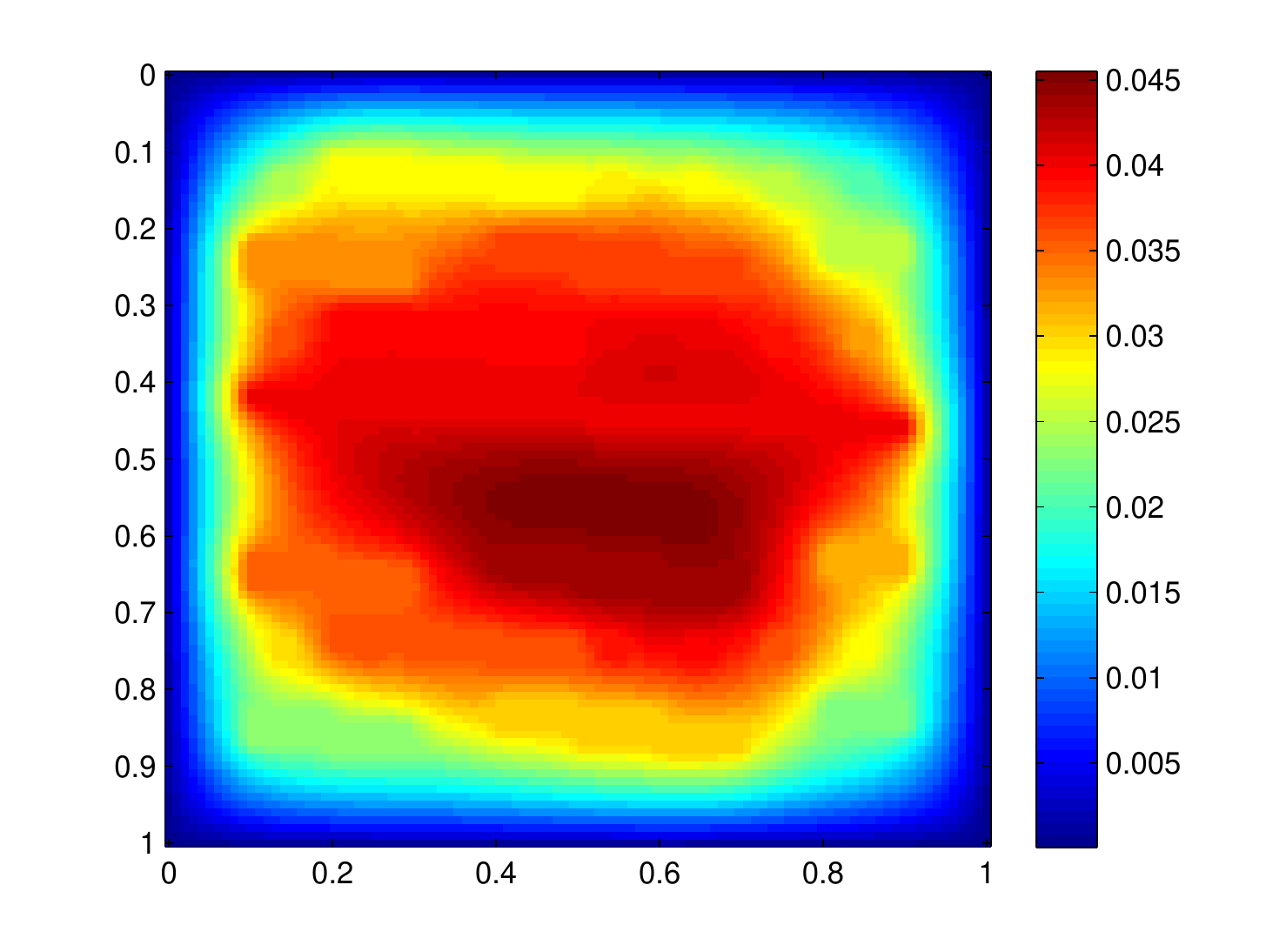}}
	\caption{ Comparison of the coarse-scale solutions with the reference (fine-scale) solution, $N=10$, $\eta=10^4$, model 1.}
	\label{model1_pressure} 
\end{figure}

\begin{figure}[H]
	\centering
	\subfigure[ Fine-scale solution]{
		\includegraphics[width=3.0in]{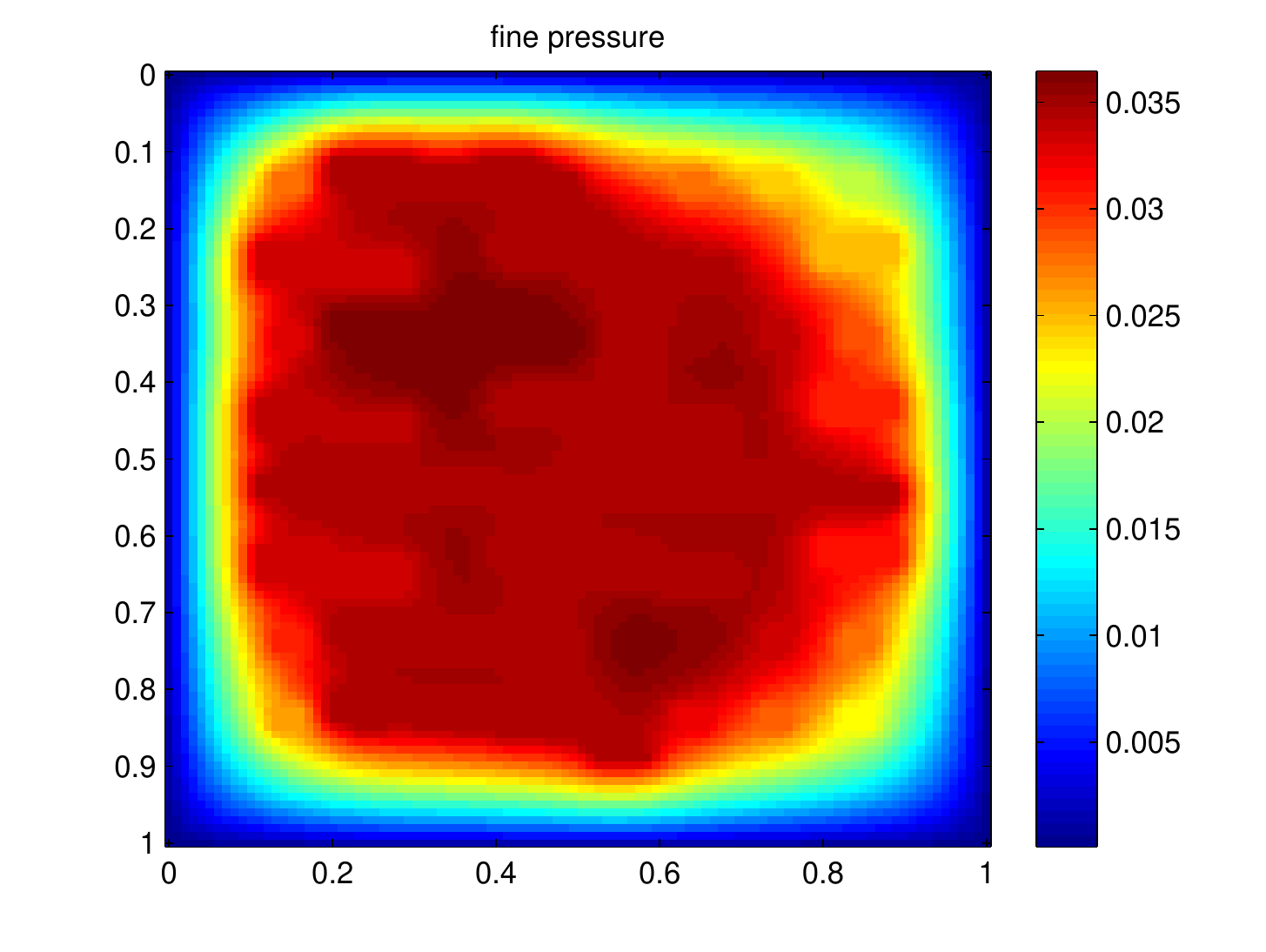}}
	\subfigure[Coarse-scale solution(1 basis)]{
		\includegraphics[width=3.0in]{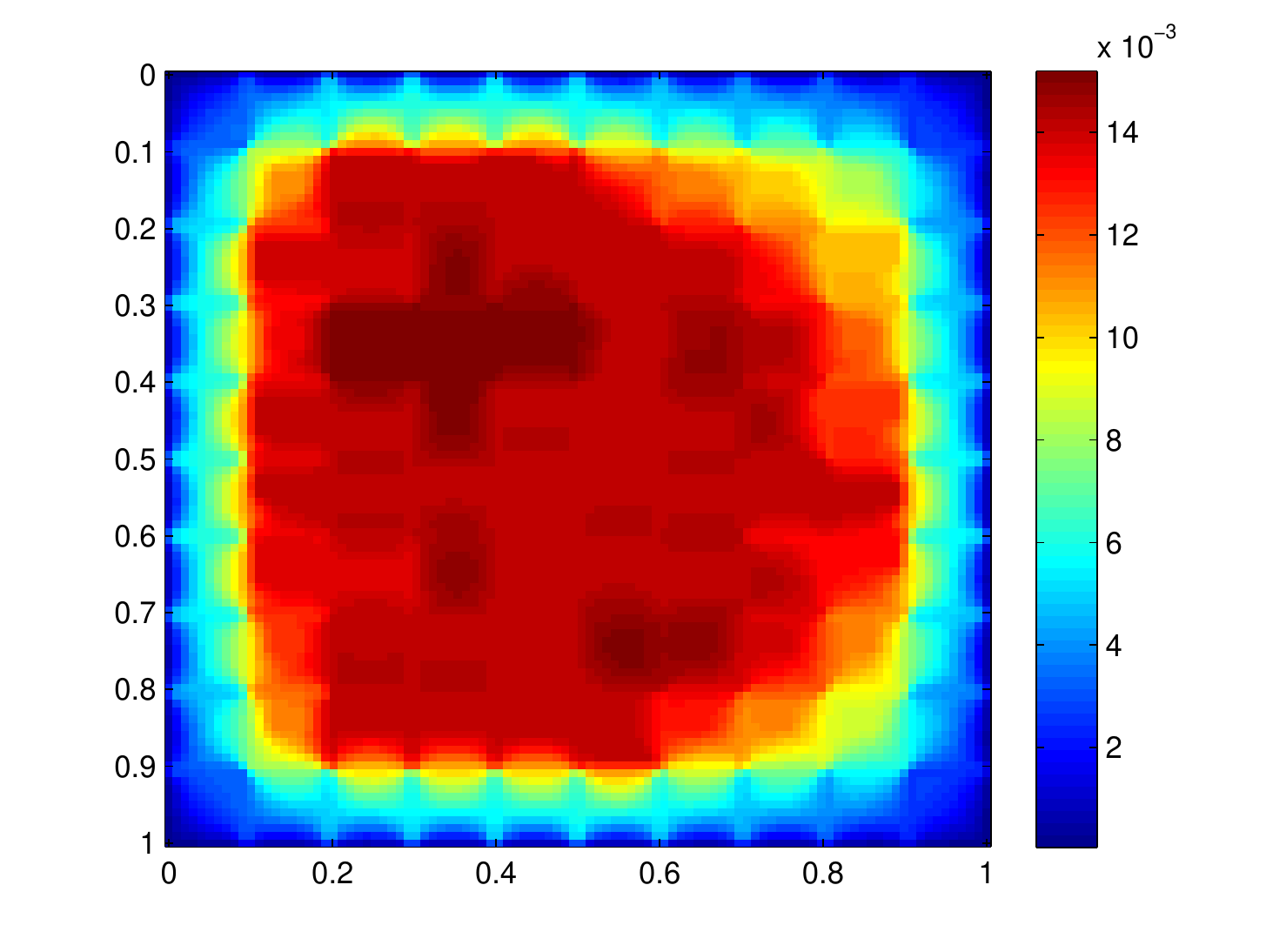}}
	\subfigure[Coarse-scale solution(3 basis)]{
		\includegraphics[width=3.0in]{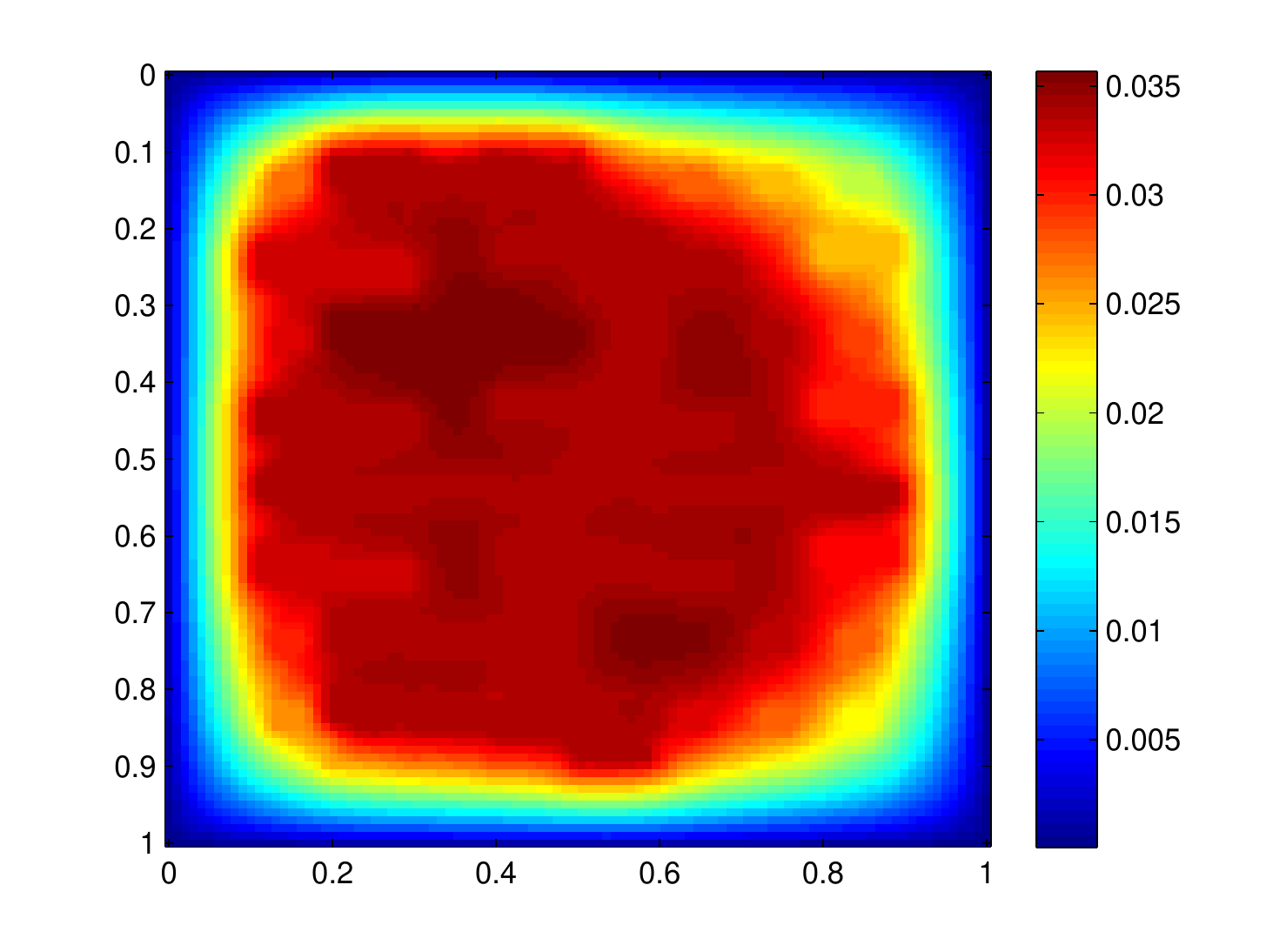}}
	\subfigure[Coarse-scale solution(5 basis)]{
		\includegraphics[width=3.0in]{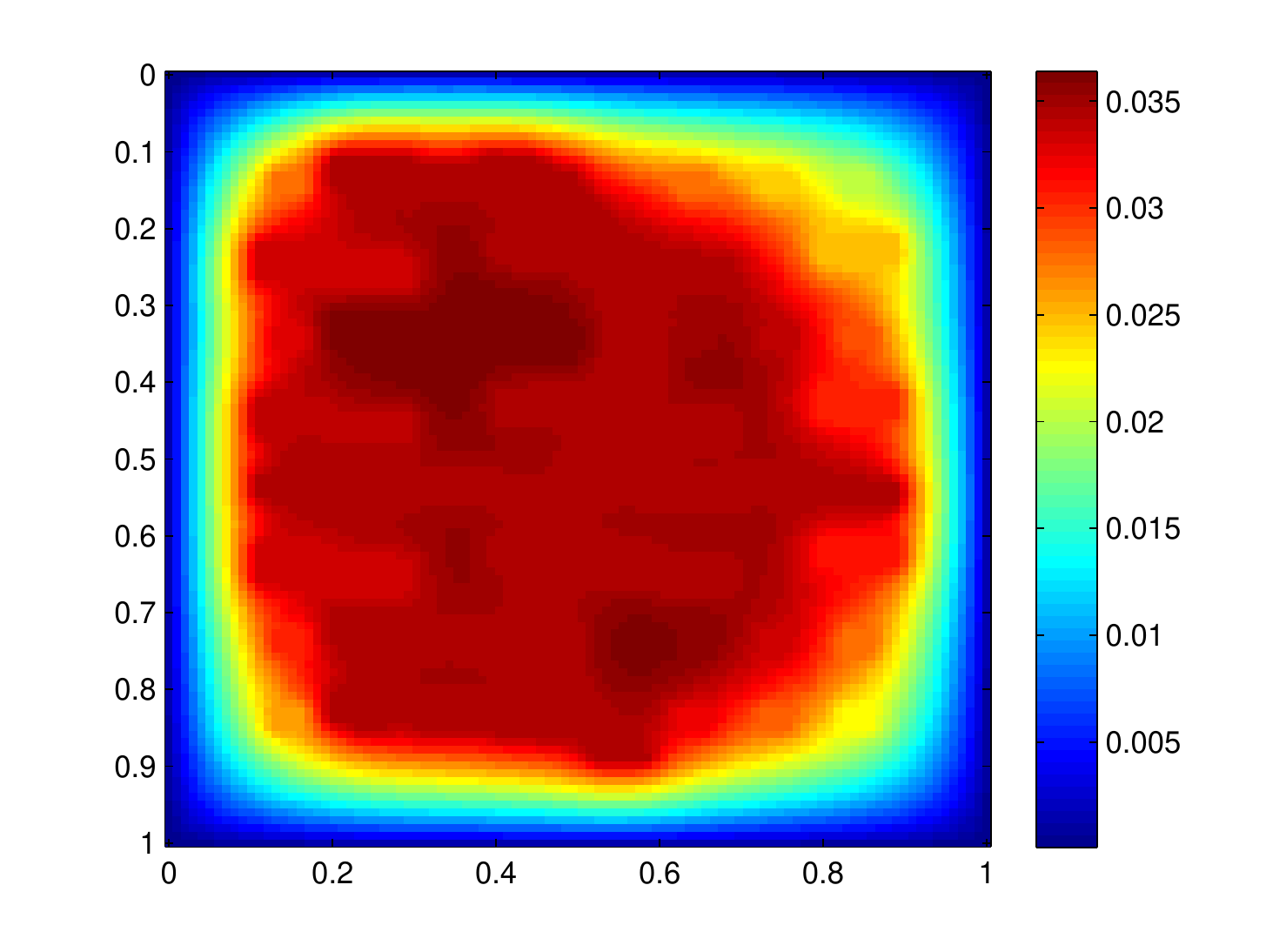}}
	\caption{ Comparison of the coarse-scale solutions with the reference (fine-scale) solution, $N=10$, $\eta=10^4$, model 2.}
	\label{model2_pressure} 
\end{figure}

\begin{figure}[H]
  \begin{center}
   \includegraphics[width=5in]{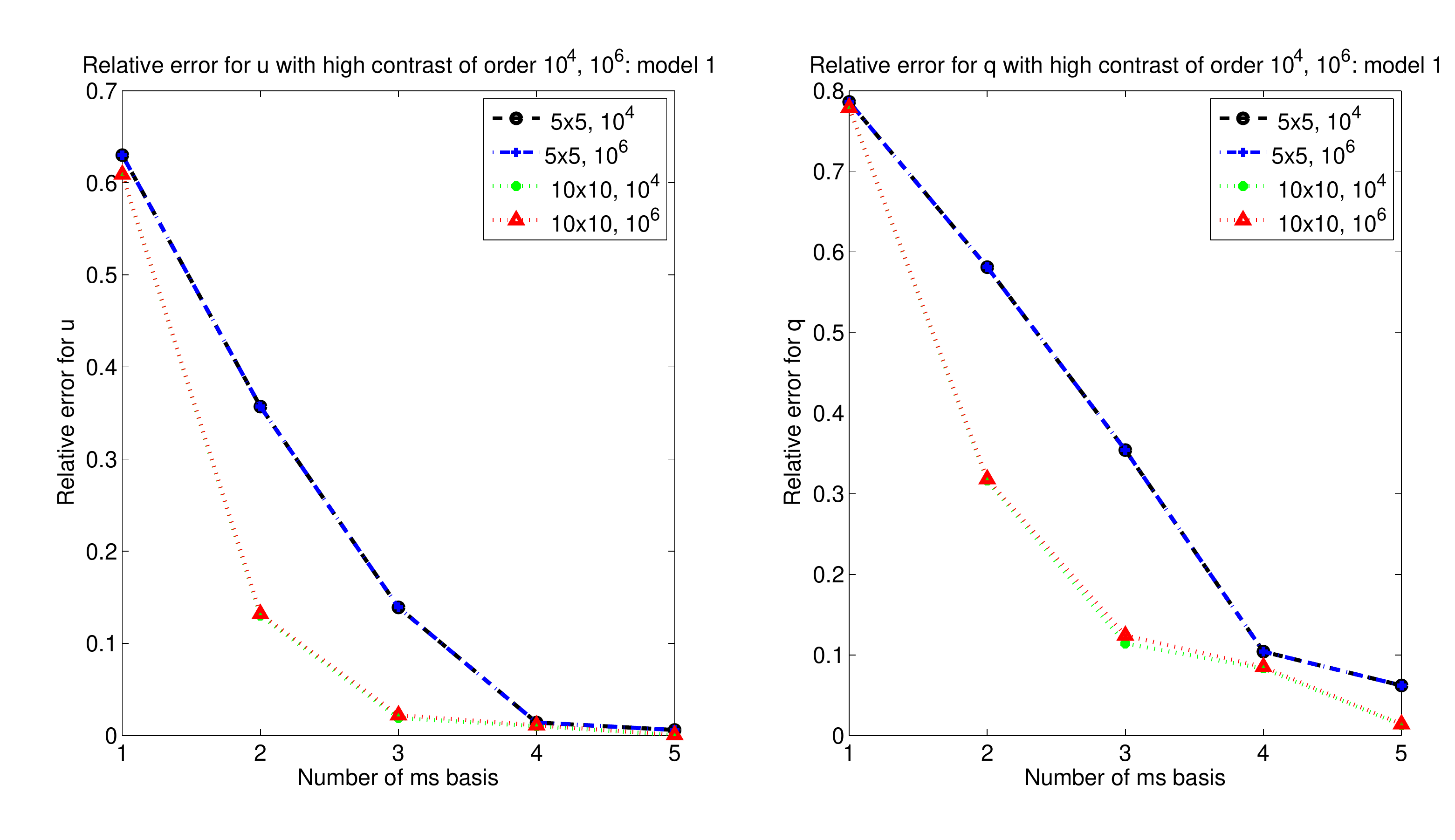}
   \end{center}
   \caption{Relative error for $u$ (left), $q$ (right) with contrast order $\eta=10^4 $ and $\eta=10^6$ for model 1, 
   	basis generation case 2.}
   \label{contrast_model1}
\end{figure}

\begin{figure}[H]
  \begin{center}
   \includegraphics[width=5in]{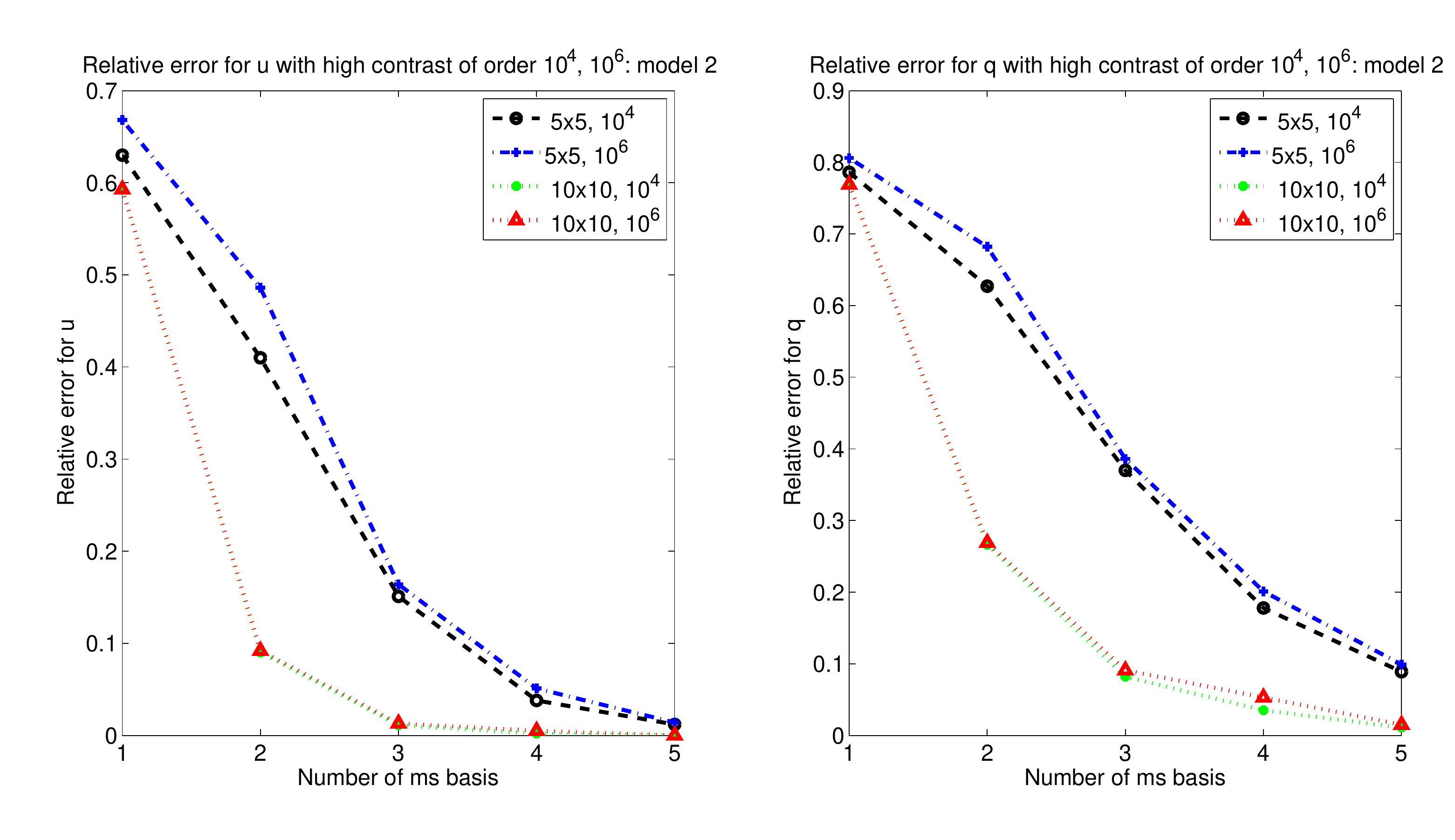}
   \end{center}
   \caption{Relative error for $u$ (left), $q$ (right) with contrast order $\eta=10^4 $ and $\eta=10^6$ for model 2, 
   	basis generation case 2.}
   \label{contrast_model2}
\end{figure}

\subsection{Preconditioner}
\par In this subsection, we present the numerical results of using 
multiscale basis to form the coarse space for the two level additive
Schwarz domain decomposition preconditioner. We use PCG as outer accelerator if the two-level preconditioner is symmetric, otherwise
GMRES with restarted number of 2 (GMRES(2)) is applied. We adopt the techniques in \cite{ganis2009implementation} to implement the local preconditioner $M_{loc}^{-1}$. Then the dominant computation can be done offline and it can be parallelized with coloring techniques, the dimension of  each local
preconditioner is equal to the number of fine scale edges in the adjacent coarse elements, which is quite small and thus can be precomputed 
and saved. Direct solver is used to implement $M_0^{-1}$. We consider both additive and hybrid preconditioners. We are particularly interested in the robustness of the method, and in each simulation we consider three types of 
contrast to test the robustness. The initial guess is zero, and the stopping
criterion that the residual is reduced by a factor of $10^7$ in $L^2$ norm.

\par The first test involves comparing the PCG iteration number of using multiscale basis and polynomial basis to form coarse space of the coarse preconditioner, which is 
shown in Table \ref{table:pcg}. The first column gives the order of high contrast. The second column is the number of PCG iterations of using polynomial to form the coarse space. The rest of the columns give the  number of PCG iterations of using multiscale basis from Case 3 and Case 4 to form the coarse space. We can see clearly that for each case, the iteration
number of multiscale basis is generally much smaller than that of polynomial basis. If polynomial basis is used, the preconditioner is not robust with respect to 
the contrast of the media. For example, for the hybrid method in the second column, PCG iteration number is 13 for contrast $10^2$, and increases to 39 for contrast  $10^6$. While for the multiscale basis,  iteration number is almost independent of contrast. Hybrid preconditioner performs  
better than additive preconditioner in terms of iteration number, however, hybrid preconditioner requires to apply the coarse preconditioner 
twice in each iteration.
\par Next, we focus  our study on the influence of the domain size on local preconditioner.  In Table \ref{table:gmres_3}, the
GMRES(2) iteration numbers of using multiscale basis from case 3 with four types of  computational domain are presented. Two multiscale basis on each coarse edge is used to form the coarse space. The first column gives the order of high contrast. The rest of the columns give the  number of GMRES(2) outer iterations for the   four types of  computational domain.  As it is shown, a restrictive local preconditioner (cases of Domain 2, 3, 4) can reduce the iteration
number approximately by half by comparing the results of Domain 1 and the rest of the 3 domain cases. For example, the number of iteration for the additive case from Domain 1 for contrast of order  $10^2, 10^4, 10^6$ is 20, 18, 18 respectively, while from Domain 2 the number of iteration for the additive case is 10, 11, 11.  Moreover, iteration number is almost independent of contrast for all types of domains. In Table \ref{table:gmres_4}, the
GMRES(2) iteration numbers of using multiscale basis from case 4 with four types of  computational domain are presented. We observed similar results as in Table \ref{table:gmres_3}. Since  multiscale basis from case 4 uses randomized technique for snapshot, it is acceptable that the number of iterations are slightly larger than  the number of iterations of corresponding cases in in Table \ref{table:gmres_3}.
\begin{table}[H]
	\centering
	\begin{tabular}{|c|c|c|c|c|c|c|}
		\hline 
		\multirow{2}{*}{Contrast} & \multicolumn{2}{c|}{\specialcell{Polynomial}}& \multicolumn{2}{c|}{\specialcell{Case 3}} & \multicolumn{2}{c|}{ \specialcell{Case 4}} \tabularnewline
		\cline{2-7} 
	& additive &hybrid& additive &hybrid & additive &hybrid \tabularnewline
		\hline
		$10^2$&23&13&25&16&26& 16 \tabularnewline
		\hline
		$10^4$ &$>40$&28&27&16&29&17  \tabularnewline
		\hline
		$10^6$&$>40$&39&27&16&28&16 \tabularnewline
		\hline
	\end{tabular}
	\caption{PCG iteration number with different types of coarse space, 2 basis is used for coarse space, $N=5$.}
	\label{table:pcg}
\end{table}

\begin{table}[h]
	\centering
	\begin{tabular}{|c|c|c|c|c|c|c|c|c|}
		\hline 
		\multirow{2}{*}{Contrast} & \multicolumn{2}{c|}{\specialcell{Domain 1}}& \multicolumn{2}{c|}{\specialcell{Domain 2}} & \multicolumn{2}{c|}{ \specialcell{Doamin 3}}& \multicolumn{2}{c|}{ \specialcell{Domain 4}} \tabularnewline
		\cline{2-9} 
		& additive &hybrid& additive &hybrid & additive &hybrid & additive &hybrid \tabularnewline
		\hline
		$10^2$&20&10&10&5&10&5&11&6 \tabularnewline
		\hline
		$10^4$ &18&10&11&5&12&5&13&7  \tabularnewline
		\hline
		$10^6$&18&10&11&5&12&5&11&6 \tabularnewline
		\hline
	\end{tabular}
	\caption{GMRES(2) iteration number with different types of preconditioners and different local preconditioner, 2 basis is used for coarse space with multiscale basis case 3, $N=5$.}
		\label{table:gmres_3}
\end{table}

\begin{table}[H]
	\centering
	\begin{tabular}{|c|c|c|c|c|c|c|c|c|}
		\hline 
		\multirow{2}{*}{Contrast} & \multicolumn{2}{c|}{\specialcell{Domain 1}}& \multicolumn{2}{c|}{\specialcell{Domain 2}} & \multicolumn{2}{c|}{ \specialcell{Doamin 3}}& \multicolumn{2}{c|}{ \specialcell{Domain 4}} \tabularnewline
		\cline{2-9} 
		& additive &hybrid& additive &hybrid & additive &hybrid & additive &hybrid \tabularnewline
		\hline
		$10^2$&17&10&15&7&12& 8&16&9 \tabularnewline
		\hline
		$10^4$ &20&8&18&5&18&10&18&7  \tabularnewline
		\hline
		$10^6$&18&10&14&6&12&7&16&9 \tabularnewline
		\hline
	\end{tabular}
	\caption{GMRES(2) iteration number with different types of preconditioners and different local preconditioner, 2 basis is used for coarse space with multiscale basis case 4, $N=5$.}
		\label{table:gmres_4}
\end{table}

\section{Conclusion}
 We have developed  an enriched multiscale mortar mixed finite element method for elliptic problems. The method is based on a mixed formulation of the problem, the concepts of domain decomposition, and mortar techniques.  The multiscale basis functions are constructed  from local problems. This method fully resolve the problem within the subdomains and glues them together with a  coarse mortar finite element space.  Further, we design both two-level {\it additive}, {\it hybrid} preconditioners which can be used within a Krylov accelerator such as PCG or GMRES. These two-level preconditioners consist of a local smoothing preconditioner based on block Jacobi(BJ), blocked by subdomain interfaces, and a coarse preconditioner based on subdomain interfaces using the enriched multiscale mortar method. Finally we present some numerical examples to demonstrate the performance of the method. 

	\bibliographystyle{plain}   
	\bibliography{references}

\begin{thebibliography}{10}

\bibitem{hybridmixed2013}
Rodolfo Araya, Christopher Harder, Diego Paredes, and Fr{\'e}d{\'e}ric
  Valentin.
\newblock Multiscale hybrid-mixed method.
\newblock {\em SIAM Journal on Numerical Analysis}, 51(6):3505--3531, 2013.

\bibitem{Arbogast_two_scale_04}
T.~Arbogast.
\newblock Analysis of a two-scale, locally conservative subgrid upscaling for
  elliptic problems.
\newblock {\em SIAM J. Numer. Anal.}, 42(2):576--598 (electronic), 2004.

\bibitem{Arbogast_PWY_07}
T.~Arbogast, G.~Pencheva, M.F. Wheeler, and I.~Yotov.
\newblock A multiscale mortar mixed finite element method.
\newblock {\em Multiscale Model. Simul.}, 6(1):319--346, 2007.

\bibitem{arbogast2000mixed}
Todd Arbogast, Lawrence~C Cowsar, Mary~F Wheeler, and Ivan Yotov.
\newblock Mixed finite element methods on nonmatching multiblock grids.
\newblock {\em SIAM Journal on Numerical Analysis}, 37(4):1295--1315, 2000.

\bibitem{arbogast2013ms}
Todd Arbogast and Hailong Xiao.
\newblock A multiscale mortar mixed space based on homogenization for
  heterogeneous elliptic problems.
\newblock {\em SIAM Journal on Numerical Analysis}, 51(1):377--399, 2013.

\bibitem{arbogast2015two}
Todd Arbogast and Hailong Xiao.
\newblock Two-level mortar domain decomposition preconditioners for
  heterogeneous elliptic problems.
\newblock {\em Computer Methods in Applied Mechanics and Engineering},
  292:221--242, 2015.

\bibitem{cai1999restricted}
Xiao-Chuan Cai and Marcus Sarkis.
\newblock A restricted additive {S}chwarz preconditioner for general sparse
  linear systems.
\newblock {\em Siam journal on scientific computing}, 21(2):792--797, 1999.

\bibitem{calo2014randomized}
V.M. Calo, Y.~Efendiev, J.~Galvis, and G~Li.
\newblock G. {R}andomized oversampling for generalized multiscale finite
  element methods.
\newblock {\em Multiscale Model. Simul.}, 14.

\bibitem{chen2003mixed}
Zhiming Chen and Thomas Hou.
\newblock A mixed multiscale finite element method for elliptic problems with
  oscillating coefficients.
\newblock {\em Mathematics of Computation}, 72(242):541--576, 2003.

\bibitem{chung2016adaptive}
Eric Chung, Yalchin Efendiev, and Thomas~Y Hou.
\newblock Adaptive multiscale model reduction with generalized multiscale
  finite element methods.
\newblock {\em Journal of Computational Physics}, 320:69--95, 2016.

\bibitem{chung2015mixed}
Eric~T Chung, Yalchin Efendiev, and Chak~Shing Lee.
\newblock Mixed generalized multiscale finite element methods and applications.
\newblock {\em Multiscale Modeling \& Simulation}, 13(1):338--366, 2015.

\bibitem{chung2015perforated}
Eric~T Chung, Yalchin Efendiev, Guanglian Li, and Maria Vasilyeva.
\newblock Generalized multiscale finite element methods for problems in
  perforated heterogeneous domains.
\newblock {\em Applicable Analysis}, pages 1--26, 2015.

\bibitem{chung2013sub}
Eric~T Chung and Wing~Tat Leung.
\newblock A sub-grid structure enhanced discontinuous {G}alerkin method for
  multiscale diffusion and convection-diffusion problems.
\newblock {\em Communications in Computational Physics}, 14(02):370--392, 2013.

\bibitem{chung2016mixed}
Eric~T Chung, Wing~Tat Leung, and Maria Vasilyeva.
\newblock Mixed {GM}s{FEM} for second order elliptic problem in perforated
  domains.
\newblock {\em Journal of Computational and Applied Mathematics}, 304:84--99,
  2016.

\bibitem{hdgunified}
Bernardo Cockburn, Jayadeep Gopalakrishnan, and Raytcho Lazarov.
\newblock Unified hybridization of discontinuous {G}alerkin, mixed, and
  continuous {G}alerkin methods for second order elliptic problems.
\newblock {\em SIAM Journal on Numerical Analysis}, 47(2):1319--1365, 2009.

\bibitem{dolean2012analysis}
Victorita Dolean, Fr{\'e}d{\'e}ric Nataf, Robert Scheichl, and Nicole Spillane.
\newblock Analysis of a two-level {S}chwarz method with coarse spaces based on
  local {D}irichlet-to-{N}eumann maps.
\newblock {\em Computational Methods in Applied Mathematics Comput. Methods
  Appl. Math.}, 12(4):391--414, 2012.

\bibitem{dur91}
L.J. Durlofsky.
\newblock Numerical calculation of equivalent grid block permeability tensors
  for heterogeneous porous media.
\newblock {\em Water Resour. Res.}, 27:699--708, 1991.

\bibitem{egw10}
Y.~Efendiev, J.~Galvis, and X.H. Wu.
\newblock Multiscale finite element methods for high-contrast problems using
  local spectral basis functions.
\newblock {\em Journal of Computational Physics}, 230:937--955, 2011.

\bibitem{efendiev2013generalized}
Yalchin Efendiev, Juan Galvis, and Thomas~Y Hou.
\newblock Generalized multiscale finite element methods ({GM}s{FEM}).
\newblock {\em Journal of Computational Physics}, 251:116--135, 2013.

\bibitem{efendiev2012robust}
Yalchin Efendiev, Juan Galvis, Raytcho Lazarov, and Joerg Willems.
\newblock Robust domain decomposition preconditioners for abstract symmetric
  positive definite bilinear forms.
\newblock {\em ESAIM: Mathematical Modelling and Numerical Analysis},
  46(5):1175--1199, 2012.

\bibitem{efendiev2009multiscale}
Yalchin Efendiev and Thomas~Y Hou.
\newblock {\em Multiscale finite element methods: theory and applications},
  volume~4.
\newblock Springer Science \& Business Media, 2009.

\bibitem{ms_hdg2}
Yalchin Efendiev, Raytcho Lazarov, Minam Moon, and Ke~Shi.
\newblock A spectral multiscale hybridizable discontinuous {G}alerkin method
  for second order elliptic problems.
\newblock {\em Computer Methods in Applied Mechanics and Engineering},
  292:243--256, 2015.

\bibitem{ms_hdg1}
Yalchin Efendiev, Raytcho Lazarov, and Ke~Shi.
\newblock A multiscale {HDG} method for second order elliptic equations. {P}art
  {I}. {P}olynomial and homogenization-based multiscale spaces.
\newblock {\em SIAM Journal on Numerical Analysis}, 53(1):342--369, 2015.

\bibitem{galvis2010domain}
Juan Galvis and Yalchin Efendiev.
\newblock Domain decomposition preconditioners for multiscale flows in
  high-contrast media.
\newblock {\em SIAM J. Multiscale Model. Simulg}, 8(4):1461--1483, 2010.

\bibitem{galvis2010domain2}
Juan Galvis and Yalchin Efendiev.
\newblock Domain decomposition preconditioners for multiscale flows in high
  contrast media: reduced dimension coarse spaces.
\newblock {\em SIAM J. Multiscale Model. Simulg}, 8(5):1621--1644, 2010.

\bibitem{ganis2009implementation}
Benjamin Ganis and Ivan Yotov.
\newblock Implementation of a mortar mixed finite element method using a
  multiscale flux basis.
\newblock {\em Computer Methods in Applied Mechanics and Engineering},
  198(49):3989--3998, 2009.

\bibitem{gao2016application}
Longfei Gao, Xiaosi Tan, and Eric~T Chung.
\newblock Application of the generalized multiscale finite element method in
  parameter-dependent {PDE} simulations with a variable-separation technique.
\newblock {\em Journal of Computational and Applied Mathematics}, 300:183--191,
  2016.

\bibitem{graham2007dd}
IG~Graham, PO~Lechner, and Robert Scheichl.
\newblock Domain decomposition for multiscale {PDE}s.
\newblock {\em Numerische Mathematik}, 106(4):589--626, 2007.

\bibitem{hybridmixed2015}
Christopher Harder, Diego Paredes, and Fr{\'e}d{\'e}ric Valentin.
\newblock A family of multiscale hybrid-mixed finite element methods for the
  {D}arcy equation with rough coefficients.
\newblock {\em Journal of Computational Physics}, 245:107--130, 2013.

\bibitem{hughes98}
T.~Hughes, G.~Feijoo, L.~Mazzei, and J.~Quincy.
\newblock The variational multiscale method - a paradigm for computational
  mechanics.
\newblock {\em Comput. Methods Appl. Mech. Engrg.}, 166:3--24, 1998.

\bibitem{ref:Iliev2011}
O.~Iliev, R.~Lazarov, and J.~Willems.
\newblock Variational multiscale finite element method for flows in highly
  porous media.
\newblock {\em Multiscale Modeling \& Simulation}, 9.4:1350--1372, 2011.

\bibitem{jennylt03}
P.~Jenny, S.H. Lee, and H.~Tchelepi.
\newblock Multi-scale finite volume method for elliptic problems in subsurface
  flow simulation.
\newblock {\em J. Comput. Phys.}, 187:47--67, 2003.

\bibitem{kim2016bddc}
Hyea~Hyun Kim, Eric Chung, and Junxian Wang.
\newblock {BDDC} and {FETI}-{DP} algorithms with adaptive coarse spaces for
  three-dimensional elliptic problems with oscillatory and high contrast
  coefficients.
\newblock {\em arXiv preprint arXiv:1606.07560}, 2016.

\bibitem{kim2016sdg}
Hyea~Hyun Kim, Eric Chung, and Chenxiao Xu.
\newblock A {BDDC} algorithm with adaptive primal constraints for staggered
  discontinuous {G}alerkin approximation of elliptic problems with highly
  oscillatory coefficients.
\newblock {\em To appear in J. Comp. Appl. Math.}

\bibitem{kim2015bddc}
Hyea~Hyun Kim and Eric~T Chung.
\newblock A {BDDC} algorithm with enriched coarse spaces for two-dimensional
  elliptic problems with oscillatory and high contrast coefficients.
\newblock {\em SIAM J. Multiscale Model. Simulg}, 13(2):571--593, 2015.

\bibitem{klawonn2015feti}
Axel Klawonn, Patrick Radtke, and Oliver Rheinbach.
\newblock {FETI}-{DP} methods with an adaptive coarse space.
\newblock {\em SIAM Journal on Numerical Analysis}, 53(1):297--320, 2015.

\bibitem{mu2016weak}
Lin Mu, Junping Wang, and Xiu Ye.
\newblock A weak {G}alerkin generalized multiscale finite element method.
\newblock {\em Journal of Computational and Applied Mathematics}, 305:68--81,
  2016.

\bibitem{oh2013overlapping}
Duk-Soon Oh.
\newblock An overlapping {S}chwarz algorithm for {R}aviart-{T}homas vector
  fields with discontinuous coefficients.
\newblock {\em SIAM Journal on Numerical Analysis}, 51(1):297--321, 2013.

\bibitem{oh2016bddc}
Duk-Soon Oh, Olof~B Widlund, Stefano Zampini, and Clark~R Dohrmann.
\newblock {BDDC} algorithms with deluxe scaling and adaptive selection of
  primal constraints for {R}aviart-{T}homas vector fields.
\newblock Technical report, tech. rep., Courant Institute, New York University,
  2015. TR2015-978, 2016.

\bibitem{balancing2003}
Gergina Pencheva and Ivan Yotov.
\newblock Balancing domain decomposition for mortar mixed finite element
  methods.
\newblock {\em Numerical linear algebra with applications}, 10(1-2):159--180,
  2003.

\bibitem{sarkis1997dd}
Marcus Sarkis.
\newblock Nonstandard coarse spaces and {S}chwarz methods for elliptic problems
  with discontinuous coefficients using non-conforming elements.
\newblock {\em Numerische Mathematik}, 77(3):383--406, 1997.

\bibitem{spillane2014abstract}
Nicole Spillane, Victorita Dolean, Patrice Hauret, Fr{\'e}d{\'e}ric Nataf,
  Clemens Pechstein, and Robert Scheichl.
\newblock Abstract robust coarse spaces for systems of {PDE}s via generalized
  eigenproblems in the overlaps.
\newblock {\em Numerische Mathematik}, 126(4):741--770, 2014.

\bibitem{toselli2005domain}
Andrea Toselli and Olof Widlund.
\newblock {\em Domain decomposition methods: algorithms and theory}, volume~34.
\newblock Springer, 2005.

\bibitem{Wheeler_MS_mortar_11}
M.F. Wheeler, T.~Wildey, and I.~Yotov.
\newblock A multiscale preconditioner for stochastic mortar mixed finite
  elements.
\newblock {\em Comput. Methods Appl. Mech. Engrg.}, 200(9-12):1251--1262, 2011.

\bibitem{Wheeler_mortar_MS_12}
M.F. Wheeler, G.~Xue, and I.~Yotov.
\newblock A multiscale mortar multipoint flux mixed finite element method.
\newblock {\em ESAIM Math. Model. Numer. Anal.}, 46(4):759--796, 2012.

\bibitem{weh02}
X.H. Wu, Y.~Efendiev, and T.Y. Hou.
\newblock Analysis of upscaling absolute permeability.
\newblock {\em Discrete and Continuous Dynamical Systems, Series B.},
  2:158--204, 2002.

\bibitem{xiao2013multiscale}
Hailong Xiao.
\newblock {\em Multiscale mortar mixed finite element methods for flow problems
  in highly heterogeneous porous media}.
\newblock PhD thesis, 2013.

\end{thebibliography}

\end{document}